%% file: manuscript.tex
\documentclass[preprint,12pt, authoryear]{elsarticle}
\pdfoutput=1
\usepackage{amssymb}
\usepackage{amsmath}
\usepackage{amsthm}
\usepackage{hyperref}
\usepackage{subfigure}% subcaptions for subfigures
\usepackage{subfigmat}% matrices of similar subfigures, aka small multiples
\usepackage{amsmath,amsfonts}
\newtheorem{assumption}{Assumption}
\newtheorem{example}{Example}
\input{pauls-macros}

\journal{CMAME}

\begin{document}

\begin{frontmatter}

\title{A density-matching approach for optimization under uncertainty}

\author{Pranay Seshadri}
\ead{ps583@cam.ac.uk}
\author{Paul Constantine}
\author{Gianluca Iaccarino}
\author{Geoffrey Parks}
\address[gp]{University of Cambridge, Cambridge CB2 1PZ, U.K.}
\address[pc]{Colorado School of Mines, Golden CO 80401, U.S.A.}
\address[ps]{Stanford University, Stanford, CA 94043, U.S.A.}

\begin{abstract}
Modern computers enable methods for design optimization that account for uncertainty in the system---so-called \textit{optimization under uncertainty}. We propose a metric for OUU that measures the distance between a designer-specified probability density function of the system response (the \textit{target}) and system response's density function at a given design. We study an OUU formulation that minimizes this distance metric over all designs. We discretize the objective function with numerical quadrature and approximate the response density function with a Gaussian kernel density estimate. We offer heuristics for addressing issues that arise in this formulation, and we apply the approach to a CFD-based airfoil shape optimization problem. We qualitatively compare the density-matching approach to a multi-objective robust design optimization to gain insight into the method. 
\end{abstract}


\end{frontmatter}

\input{introduction}

\input{math}

\input{discretization}

\input{heuristics}

\input{examples}

\input{conclusions}

\section*{Acknowledgments}
This research was funded through a Dorothy Hodgkin Postgraduate Award, which is jointly sponsored by the Engineering and Physical Sciences Research Council (EPSRC) (UK) and Rolls-Royce plc. The first author would like to acknowledge the financial assistance provided by the Center for Turbulence Research at Stanford University and St.~Edmund's College, Cambridge. The authors would like to thank Shahrokh Shahpar of Rolls-Royce plc for his advice on various aspects of this work. The authors also thank the reviewers for their suggestions and comments, which improved the overall quality of this manuscript. The second author's work is supported by the U.S. Department of Energy Office of Science, Office of Advanced Scientific Computing Research, Applied Mathematics program under Award Number DE-SC-0011077.

%%\section*{References}
\bibliographystyle{elsarticle-harv}
\bibliography{references}

\end{document}

%% file: pauls-macros.tex
\usepackage{bm}

\graphicspath{{./}{./figs/}}

  \def\clap#1{\hbox to 0pt{\hss#1\hss}}

\providecommand{\mat}[1]{\bm{#1}}%
\renewcommand{\vec}[1]{\mathbf{#1}}

\providecommand{\mF}{\ensuremath{\mat{F}}}

\providecommand{\mK}{\ensuremath{\mat{K}}}

\providecommand{\mW}{\ensuremath{\mat{W}}}

\providecommand{\ve}{\ensuremath{\vec{e}}}
\providecommand{\vf}{\ensuremath{\vec{f}}}

\providecommand{\vq}{\ensuremath{\vec{q}}}

\providecommand{\vt}{\ensuremath{\vec{t}}}

\newcommand{\sS}{\mathcal{S}}

\newcommand{\bmat}[1]{\begin{bmatrix}#1\end{bmatrix}}

\newcommand{\argmin}[1]{\underset{#1}{\mathrm{argmin}}}

%% file: introduction.tex
\section{Introduction}
\label{sec:intro}

\noindent Modern computing power opens the possibility for industrial-scale design optimization with high-fidelity numerical simulations of physical systems. Simulation-based design is found in aircraft~\cite{Airbus}, engine~\cite{Shahrokh}, automotive~\cite{auto} and shipping~\cite{Shipping} industries, among many others. To optimize, designers must precisely specify operating scenarios and manufactured production. Off-design operation and manufacturing tolerances are typically incorporated afterward. A more complete perspective on design optimization accounts for these uncertainties, e.g., by employing statistical performance metrics within the design optimization. This perspective leads to \textit{optimization under uncertainty} (OUU).

The computational engineering literature is chock full of formulations and approaches for OUU. Allen and Maute~\cite{Allen} give an excellent overview that broadly categorizes these formulations as either \textit{robust design optimization} (RBO) or \textit{reliability-based design optimization} (RBDO). The essential idea behind RBO formulations is to simultaneously maximize a statistical measure of the system performance (e.g., the mean) while minimizing a statistical measure of system variability (e.g., the variance), thus improving robustness to variability in operating conditions. The optimization is often formulated with multiple objective functions (e.g., maximize mean and minimize variance), which leads to a Pareto front of solutions representing a trade-off between robustness and performance. Alternative formulations treat performance as the objective function and robustness as a constraint or vice versa. Some applications of RBO include the design of Formula One brake ducts~\cite{Axerio}, compressor blades~\cite{Seshadri_ASME}, compression systems~\cite{Ghisu}, airfoils~\cite{Gianluca_windmill}, and structures~\cite{Dolt}. The RBDO formulations seek designs that satisfy reliability criteria, such as maintaining a sufficiently small probability of failure, while minimizing a cost function of the design~\cite{Frangopol}. Estimating the failure probabilities within the optimization with randomized methods (e.g., Monte Carlo) can be prohibitively expensive for large-scale models; several methods exist for approximating regions of low failure probability~\cite{Bichon}. Engineering examples of RBDO include transonic compressors~\cite{Lian}, aeroelasticity~\cite{Missoum}, structures~\cite{Allen}, and vehicle crash worthiness~\cite{Carcrash}.

The statistical measures in the RDO and RBDO objective functions and constraints are typically low-order moments---e.g., mean and variance---or probabilities associated with the system response. The chosen statistical measures affect the optimal design, so they must be chosen carefully for each specific application. 

In this paper, we propose an alternative statistical measure that can be used in the RDO and RBDO formulations. We assume the designer has described the desired system performance as a full probability density function (pdf), which we call the \textit{target pdf}, and we seek to minimize the distance between the design-dependent response pdf and the target pdf. Mathematically, we treat the target pdf as given; it is not a tunable parameter. In any real-world scenario, this pass-the-buck attitude places tremendous responsibility on the designer to devise the perfect target pdf. We expect that a practical methodology including the proposed statistical measure will involve some back and forth between designer and optimizer to devise the most appropriate target pdf. Using a designer-specified response pdf has some precedent in the OUU literature. Rangavajhala and Mahadevan~\cite{Rang} assume a designer-specified pdf in their \textit{optimum threshold design}, which finds thresholds that satisfy the given joint probability while allowing for preferences among multiple objectives.

We present a single-objective OUU formulation where the distance between target and response pdfs is the objective function. We explore some interesting properties of this optimization problem, namely how the objective's gradient behaves when the two pdfs are not sufficiently large on the same support (section \ref{sec:math}). We propose a convergent discretization of the objective function---based on numerical quadrature and kernel density estimation---that produces a continuous approximation well-suited for gradient-based optimization (section \ref{sec:discretization}). Our prior work uses histograms to approximate the response pdf, which leads to a less scalable optimization problem with integer variables~\cite{Aggressive}. There are some drawbacks to the density-matching formulation, and offer heuristics for addressing these drawbacks in section \ref{sec:heuristics}. In section \ref{sec:examples}, we test the formulation on an algebraic test problem and a shape optimization problem with the NACA0012 airfoil. In the latter case, we qualitatively compare the optimal designs to those generated by a multi-objective RDO strategy.

%% file: math.tex
\section{Mathematical formulation}
\label{sec:math}
\noindent Consider a function $f=f(s,\omega)$ that represents the response of a physical model with design variables $s\in\mathcal{S}\subseteq\mathbb{R}^{n}$ and random variables $\omega\in\Omega\subseteq\mathbb{R}^m$; the random variables represent the uncertainty in the physical system. The space $\sS$ encodes the application-specific constraints on the design variables, e.g., bounds or linear inequality constraints. We assume that $\omega$ are defined on a probability space with sample space $\Omega$ and probability density function $p=p(\omega)$, which encode all available knowledge about the system's uncertainties.\footnote{The final results depend on $\Omega$ and $p(\omega)$. If multiple probability density functions are consistent with the available information, then one should check the sensitivity of the results to perturbations in these quantities.} We assume that $f$ is scalar-valued, $f\in\mathcal{F}\subseteq\mathbb{R}$, though this can be generalized. We also assume that $f$ is continuous in both $s$ and $\omega$.  For a fixed  $s\in\mathcal{S}$, let $q_{s}:\mathcal{F}\rightarrow\mathbb{R}_{+}$ be a probability density function of $f(s,\omega)$; we assume that $f(s,\omega)$ admits a square-integrable pdf for all values $s$ in the design space $\mathcal{S}$. The shape of $q_s$ will be different for different values of $s$. 

The given target pdf expresses the designer's desired system performance accounting for uncertainty in operating conditions. Denote the target pdf by $t:\mathbb{R}\rightarrow\mathbb{R}_+$, which we assume is square-integrable. To find the values of the design variables $s$ that bring the system's response as close as possible to the designer's target, we pose the following optimization problem:
\begin{equation}
\label{eq:opt}
s^\ast \;=\; \argmin{s\in\sS} \; d(t,q_s),
\end{equation}
where $d(\cdot,\cdot)$ is a distance metric between two comparable probability density functions. The values $s^\ast$ correspond to the optimal design under uncertainty.

A few comments on this optimization problem are in order. Since $d$ is a distance metric, $d\geq 0$. However, $d(t,q_s)$ is not generally a convex function of $s$. Therefore, $s^\ast$ may not be unique, and the optimization problem may need a regularization term to make it well-posed (e.g., Tikhonov regularization). 

The minimum value of the objective function $d(t,q_{s^\ast})$ measures how well the optimal design meets the designer's specifications. A non-zero value at the minimum means that the model can be improved, e.g., by incorporating more controls or otherwise modifying the relationship between the design variables and the system behavior. If the minimizing design is deemed too far from the target, then the designer may request a radical redesign of the system; by introducing additional design variables or increase their ranges, that allows the model to get closer to her specifications.

The formulation in \eqref{eq:opt} uses $d$ as the sole objective function, and \eqref{eq:opt} has no constraints that depend on the uncertainties. We study this formulation because of its simplicity. One could use the metric $d$ as a one objective in a multi-objective RDO formulation or as a measure of reliability in an RBDO formulation. We do not pursue these ideas in this paper. 

There are many possible choices for the distance metric $d$; Gibbs and Su~\cite{Gibbs2002} review several metrics and the relationships between them. To enable efficient, scalable gradient-based methods for the optimization \eqref{eq:opt}, we choose the differentiable squared $L_2$-norm, 
\begin{equation}
\label{eq:obj}
d(t,q_s) \;=\; \int_{-\infty}^{\infty} (t(f) - q_s(f))^2\,df.
\end{equation}
This integral is finite by the square-integrability assumption on $t$ and $q_s$. 

\subsection{The trouble with non-overlapping response and target pdfs}
\label{sec:overlap}
\noindent Something peculiar happens to $d$ from \eqref{eq:obj} when the supports of $t$ and $q_s$ do not overlap---i.e., $t$ is zero if $q_s$ is positive and vice versa. Expanding the integrand in \eqref{eq:obj},
\begin{equation}
d(t,q_s) \;=\; \int_{-\infty}^{\infty} t(f)^2\,df
- 2\int_{-\infty}^{\infty} t(f)\,q_s(f)\,df
+ \int_{-\infty}^{\infty} q_s(f)^2\,df.
\end{equation}
Since the target $t$ is independent of the design variables $s$, the minimizer of $d$ is the same as the minimizer of $d'$ defined as
\begin{equation}
\label{eq:objp1}
d'(t,q_s) \;:=\; \int_{-\infty}^{\infty} q_s(f)^2\,df - 2\int_{-\infty}^{\infty} t(f)\,q_s(f)\,df
\end{equation}
If the supports of $t$ and $q_s$ do not overlap, then the second term in \eqref{eq:objp1} vanishes, and
\begin{equation}
\label{eq:objp2}
d'(t,q_s) \;=\; \int_{-\infty}^{\infty} q_s(f)^2\,df.
\end{equation}
In words, when $t$ and $q_s$ do not overlap, the objective function has no information from the target $t$. The gradient of $d'$ with respect to the design variables $s$ may point in a direction that decreases $d'$, but there is no guarantee that a step along that direction in the design space moves $q_s$ closer (by the distance metric) to the target $t$. The following example illustrates this issue.

\begin{example}
Let $f(s,\omega)=s+\omega$, where $\omega$ is a random variable distributed uniformly on $[0,1]$ and $s\in[0,2]$. The response pdf $q_s(f)$ is a uniform density function on the interval $[s,s+1]$. Let the target pdf $t(f)$ be a uniform density on the interval $[2,3]$. The minimizer of \eqref{eq:opt} is $s^\ast=2$, and $d(t,q_{s^\ast})=0$. However, for $s\in[0,1)$, the first and second derivatives of $d(t,q_s)$ with respect to $s$ are zero. Thus, a gradient-based optimization method would stop at any candidate minimizer in the interval $[0,1)$. 
\end{example}

\noindent The overlap issue forces us to make the following assumption to ensure that the minimizer from a gradient-based method applied to \eqref{eq:opt} with distance metric \eqref{eq:obj} produces a response pdf with some relationship to the given target pdf. 

\begin{assumption}
\label{bigass}
Assume that for all $s\in\sS$, the intersection of the support of the target pdf $t(f)$ and the support of the design-dependent response pdf $q_s(f)$ is non-empty. 
\end{assumption}

\noindent Assumption \ref{bigass} is sufficient but not necessary; it can be relaxed to (i) the initial design point produces a response pdf whose support overlaps the target's support and (ii) all iterates of the optimization method produce response pdfs whose supports overlap the target's support. For computation, we exploit choices in the kernel density estimates to ensure that Assumption \ref{bigass} is satisfied. This approach also suggests a heuristic to accelerate the numerical optimization; see section \ref{sec:kde}. 

One final comment on the number of components in $f$: in principle, our construction can be extended to $f$'s that return a vector of responses from the system. However, this situation necessitates (i) a joint probability density for the target, (ii) a multivariate density estimation method for the response pdf, and (iii) multivariate integration to compute the distance metric. Thus, the approach suffers the dreaded curse of dimensionality as the number of components in $f$ increases. However, it worthwhile to note that it is unlikely that designers would have more than a handful of objectives to optimizer over. For such cases, existing approaches such as RDO---with a multi-objective mean-variance form---lead to multiplication of objectives, rendering a more complex optimization problem. Our strategy on the other hand, still remains a `simple' distance minimization problem.

%% file: discretization.tex
\section{Discretization and computation}
\label{sec:discretization}
\noindent Next we turn to the computational aspects of solving the optimization problem \eqref{eq:opt} using the distance metric \eqref{eq:obj}. There are two main issues to address: (i) discretizing the integral in the distance metric and (ii) estimating the response density $q_s$. 

\subsection{Discretizing the distance metric}
\noindent To avoid issues with numerical integration on unbounded domains, we assume that $f(s,\omega)$ is bounded for all $s$ and $\omega$,
\begin{equation}
\label{eq:bnds}
f_\ell \;\leq\; f(s,\omega) \;\leq\; f_u,\qquad s\in\sS,\;\omega\in\Omega.
\end{equation}
This implies that the support of $q_s(f)$ is always finite. Such an assumption is not terribly restrictive. For a particular design point, $q_s$ may have a long tail, but any computer representation of this long tail necessarily imposes finite bounds. The bounds $f_\ell$ and $f_u$ need not be tight. But finite bounds helps us devise a practical discretization. The bounds imply
\begin{equation}
\label{eq:distance}
d(t,q_s) \;=\; \int_{f_\ell}^{f_u} (t(f) - q_s(f))^2\,df 
+ \int_{-\infty}^{f_\ell} t(f)^2\,df 
+ \int_{f_u}^{\infty} t(f)^2\,df.
\end{equation}
Since the target pdf $t$ is independent of the design variables $s$, the optimization can ignore the last two terms in \eqref{eq:distance}. 

We choose an $N$-point numerical quadrature rule on the interval $[f_\ell,f_u]$ with points $\gamma_i\in[f_\ell,f_u]$ and associated weights $w_i$ with $i=1,\dots,N$. The number $N$ of points in the integration rule can be very large ($\mathcal{O}(10^6)$) without a large computational burden. We need to evaluate the given target pdf $t$ and an estimate of the response pdf $q_s$ at each quadrature node, but this is very cheap. The discretized objective function is
\begin{equation}
\label{eq:aopt}
\hat{d}(t,q_s) 
\;=\; \sum_{i=1}^N (t(\gamma_i) - q_s(\gamma_i))^2\,w_i
\;=\; (\vt-\vq_s)^T\mW(\vt-\vq_s),
\end{equation}
where
\begin{equation}
\mW=\bmat{w_1 & & \\ & \ddots & \\ & & w_N},\quad
\vt=\bmat{t(\gamma_1)\\ \vdots \\ t(\gamma_N)},\quad
\vq_s=\bmat{q_s(\gamma_1)\\ \vdots \\ q_s(\gamma_N)}.
\end{equation}
The resolution of the points $\{\gamma_i\}$ should be fine enough to resolve both the target $t$ and the response pdf $q_s$ for all $s\in\sS$. In an extreme example of insufficient resolution, the support of $t$ may be entirely inside an interval defined by two neighboring quadrature nodes. In such a case, $t$ does not affect the discretized objective $\hat{d}$ in \eqref{eq:aopt}. If the support of $t$ or a particular $q_s$ is very small relative to the interval $[f_\ell,f_u]$, then one might consider a non-uniform distribution of quadrature nodes to properly resolve the pdfs. However, as noted, an extremely fine grid does not greatly increase the cost of evaluating the discretized objective function. So resolution---and, consequently, discretization error in the integral from \eqref{eq:obj}---is not a primary concern.

The choice of quadrature rule depends on the smoothness of the target $t$ and the response $q_s$ (or, its estimate)~\cite{Davis2007}. If these pdfs are very smooth on the interval $[f_\ell,f_u]$, then one could use high-order, interpolatory quadrature rules like Gaussian quadrature or Clenshaw-Curtis quadrature~\cite{Trefethen2008}. However, say the target is a uniform density on a small interval. Then a low-order method like the trapezoidal rule may be more appropriate. We prefer a highly resolved trapezoidal rule in general; recent analysis shows that it compares well to high-order methods for smooth functions~\cite{Trefethen2014}. 

\subsection{Estimating the response density}
\label{sec:densest}

\noindent For a fixed design point $s\in\sS$, the density $q_s$ is, in general, not a known function of $f$ and must be estimated. We propose to use a kernel density estimate for $q_s$, which has several advantages. We draw a set of $M$ points $\{\omega_j\}$ independently according to the given density $p(\omega)$ on the random variables representing uncertainty.
Define the functions
\begin{equation}
f_j(s) \;=\; f(s,\omega_j),\qquad j=1,\dots,M. 
\end{equation}
For a bandwidth parameter $h$ and a radial kernel $K(\cdot)$ that depends on $h$, we approximate $q_s$ by
\begin{equation}
\label{eq:kde}
q_s(f) \;\approx\; \hat{q}_s(f) \;=\; \frac{1}{M} \sum_{j=1}^M K(f-f_j(s)).
\end{equation}
We approximate the vector $\vq_s$ from \eqref{eq:aopt} as
\begin{equation}
\label{eq:vqs}
\vq_s \;\approx\; \hat{\vq}_s \;=\; \mK\ve,\qquad \mK\in\mathbb{R}^{N\times M},
\end{equation}
where 
\begin{equation}
\mK_{ij}=\frac{1}{M} K(\gamma_i-f_j(s)),\qquad i=1,\dots,N,\;
j=1,\dots,M,
\end{equation} 
and $\ve$ is an $M$-vector of ones. For computation, we replace $\vq_s$ by $\hat{\vq}_s$ in the approximate objective function $\hat{d}$ in \eqref{eq:aopt}. 

For a sufficiently small $h$, the asymptotic mean-squared error in the kernel density estimate decreases as the number $M$ of samples increases~\cite[Chapter 6]{Scott}. In practice, one can increase the bandwidth parameter $h$ to create smooth estimates that compensate for too few samples. Such heuristics are appropriate, since our goal is not perfect representation of the response density $q_s$. Our goal is to find a design point $s^\ast$ whose corresponding response pdf is sufficiently close to the given target pdf. Nevertheless, if the number $n$ of components in $\omega$ is large, then one might be concerned that $M$ is not large enough to represent the system response over the high-dimensional space $\Omega$, resulting in a poor approximation of $q_s$---potentially poor enough to adversely affect the optimization. This concern is valid when evaluating $f$ is computationally expensive, e.g., if the system involves complex computational fluid dynamics model, thus limiting $M$. In this case, we might construct a response surface of $f(s,\omega)$ as a function of $\omega$ to sample in place of the true system response. We discuss the benefits and drawbacks of response surfaces in section \ref{sec:respsurf}. 

We propose a Gaussian kernel for the density estimate,
\begin{equation}
\label{eq:gaussk}
K(r) \;=\; \frac{1}{\sqrt{2\pi}}\exp\left(-(r/h)^2/2\right).
\end{equation}
There are two main advantages to using a Gaussian kernel. First, it is differentiable at all points in its infinite domain, so we can take the derivative of the density estimate without worrying about non-differentiability at kernel support boundaries; many compactly supported kernels do not enjoy such an advantage. Second, the infinite support of the kernel implies that, at least mathematically, we can address the non-overlapping issue discussed in section \ref{sec:overlap}. It might seem inconsistent to use a kernel with infinite support when we assume that the support of $q_s$ is bounded according to \eqref{eq:bnds}. But approximating compactly supported densities with Gaussian kernel density estimates is common; we can control the error in the tails to keep it from heavily influencing the objective function in \eqref{eq:aopt}. In section \ref{sec:kde}, we propose a heuristic that exploits the freedom in the bandwidth parameter to help ensure that the $q_s$'s kernel estimate is sufficiently large on the target pdf's support. 

\subsection{Computing the gradient}
\label{sec:grad}

\noindent We can compute the gradient of the approximate objective $\hat{d}$ with respect to the design variables $s$. For the $k$th component of $s$, denoted $s_k$, 
\begin{equation}
\frac{\partial \hat{d}}{\partial s_k}
 = 2\,(\vt-\hat{\vq}_s)^T\,\mW\,
\left(\frac{\partial\hat{\vq}_s}{\partial s_k}\right).
\end{equation}
Note that
\begin{equation}
\frac{\partial \hat{q}_s}{\partial s_k} 
\;=\;
\frac{1}{M} \sum_{j=1}^M K'(f-f_j(s))\,\frac{\partial f_j}{\partial s_k},
\end{equation}
where $K'$ is the derivative of the kernel with respect to its argument, which is easily computed from \eqref{eq:gaussk}. To reiterate, the partial derivative $\partial f_j/\partial s_k$ is the derivative of the response $f(s,\omega)$, with $\omega=\omega_j$, with respect to the $k$th design variable $s_k$. This partial derivative is a function of the design variables. Define
\begin{equation}
\vf'_k = \bmat{\frac{\partial f_1}{\partial s_k}\\ \vdots \\ \frac{\partial f_M}{\partial s_k}},\qquad
\mK'_{ij} = \frac{1}{M} K'(\gamma_i-f_j(s)).
\end{equation}
Then we can concisely write the derivative of $\hat{d}$ from \eqref{eq:aopt} with respect to the $k$th component of $s$ as
\begin{equation}
\frac{\partial \hat{d}}{\partial s_k}
\;=\;
2\,\left(\vt-\mK\ve\right)^T\,\mW\,\mK'\,\vf'_k.
\end{equation}
Define the $M\times n$ matrix $\mF'$ by
\begin{equation}
\label{eq:grads}
\mF' \;=\; 
\bmat{
\frac{\partial f_1}{\partial s_1} & \cdots & \frac{\partial f_1}{\partial s_n}\\
\vdots & \ddots & \vdots\\
\frac{\partial f_M}{\partial s_1} & \cdots & \frac{\partial f_M}{\partial s_n}\\
}.
\end{equation}
We can write the gradient of the objective function---oriented as a row vector---as
\begin{equation}
\label{eq:equation_gradients}
\nabla_s \hat{d} \;=\; 
2\,\left(\vt-\mK\ve\right)^T\,\mW\,\mK'\,\mF'.
\end{equation}
The elements of $\mK$, $\mK'$, and $\mF'$ all depend on $s$. Recall the dimensions of the terms in \eqref{eq:equation_gradients}. The gradient vector $\nabla_s \hat{d}$ has $n$ components, which is the number of random variables describing the system's uncertainty. The vector $\vt$ has $N$ components, which is the number of quadrature nodes from \eqref{eq:aopt}; we expect this to be a very large number. The matrix $\mK$ has size $N\times M$, where $M$ is the number of randomly chosen points in $\Omega$ used to estimate the pdf $q_s$. If evaluating the response is cheap, or if the response is approximated by a response surface, then $M$ may also be very large. The vector $\ve$ of ones has length $M$. The diagonal matrix $\mW$ has $N$ nonzero elements on the diagonal. The matrix $\mK'$ has size $N\times M$, and the matrix $\mF'$ has size $M\times n$. Our numerical studies have not needed special methods to evaluate the matrix-vector products in \eqref{eq:equation_gradients}. However, with a Gaussian kernel, we could perform extremely large computations (i.e., large $M$ and $N$) with a fast Gaussian transform~\cite{Greengard1991}. 

\subsection{Interfaces and cost}
\label{sec:cost}
\noindent In terms of interfaces to the simulation code, we need to evaluate (i) $f$ given $s$ and $\omega$, and (ii) the gradient $\nabla_s f$ given $s$ and $\omega$---similar to a deterministic optimization. In this sense, the approach is \emph{non-intrusive}. If we use a gradient-based optimization algorithm, such as a sequential quadratic program~\cite[Chapter 16]{Nocedal}, then each iteration uses $M$ evaluations of $f$ and its gradient with respect to $s$.

%% file: heuristics.tex
\section{Computational heuristics}
\label{sec:heuristics}

\noindent In this section, we discuss two heuristics for the optimization in \eqref{eq:opt}. The first is an approach to the kernel bandwidth selection that alleviates the non-overlapping issue discussed in section \ref{sec:overlap}. The second is the use of response surfaces in place of the true response for expensive simulations. We end this section with a short discussion of some implementation details.  

\subsection{Bandwidth parameter and the overlap problem}
\label{sec:kde}
\noindent There is a great deal of work on the proper bandwidth choice in kernel density estimation~\cite[Chapter 6]{Scott}. In most statistical inference, the data determines the bandwidth parameter~\cite{Sheather1991}. Our goal is somewhat different. Indeed, we want a reasonable estimate of the response pdf $q_s$ at a design point $s$, where we can treat the set of scalars $\{f_j(s)\}$ as data. However, we can also use the bandwidth parameter to help ensure sufficient overlap between the estimate $\hat{q}_s$ and the target pdf $t$, thus aiding the optimization. Mathematically, by using the Gaussian kernel in \eqref{eq:gaussk} with infinite support, the kernel estimate $\hat{q}_s$ is strictly positive over the the entire support of $t$. Numerically, the value of $\hat{q}_s$ may be too small on $t$'s support to produce a useful gradient. 

If $t$'s support resides in $\hat{q}_s$'s tail, then increasing the bandwidth $h$ increases $\hat{q}_s$ over $t$'s support. Loosely speaking, we can use a large $h$ to help the response pdf find the target pdf. The large $h$ produces estimates of $q_s$ that are too smooth with large halfwidths. But the large width produces useful gradients for the optimizer to help bring the response pdf closer (in terms of the $L_2$ distance) to the target. For unimodal response pdfs, this leads to more overlap between the target and the response. Once the optimizer has found a region of the design space $\sS$ where there is sufficient overlap, we reduce $h$ to data-driven values to better estimate $q_s$. 

For a kernel estimate of a Gaussian pdf using a Gaussian kernel, the optimal bandwidth is
\begin{equation}
\label{eq:scott}
h_{\text{opt}}=\left(\frac{4}{3M}\right)^{1/5}\sigma,
\end{equation}
where $\sigma$ is the Gaussian's standard deviation. The formula \eqref{eq:scott} is known as Scott's rule \cite{Bowman,Scott}. In the first few optimization iterations, we use an initial bandwidth $h=(f_u-f_\ell)/5$. Once we are satisfied that the design point $s$ yields an estimate of $q_s$ whose support sufficiently overlaps $t$, we reduce the bandwidth to $h=h_{\text{opt}}$. The results in section \ref{sec:airfoil} use this heuristic. 

\subsection{Response surfaces}
\label{sec:respsurf}
\noindent When the simulation is expensive, the number $m$ of random variables is sufficiently small, and the response $f(s,\omega)$ is a sufficiently smooth function of $\omega$, it may be more efficient to use a response surface when approximating the pdf $q_s$. Response surfaces for approximating pdfs are common in uncertainty quantification~\cite[Chapter 13]{smith2013uncertainty}. Popular response surface constructions include polynomial approximations~\cite[Chapter 7]{Dongbin}\cite[Chapter 3]{lemaitre2010spectral} and radial basis approximations~\cite{Wendland2005}. The essential idea is, for a fixed design point $s$, evaluate $f(s,\omega)$ at a few points $\omega_k\in\Omega$ with $k=1,\dots,P$, where $P$ is smaller than the number of points needed to accurately estimate the pdf $q_s$; let $f_k(s)=f(s,\omega_k)$. Most response surface constructions are linear models of the data,
\begin{equation}
f(s,\omega) \;\approx\; \tilde{f}(s,\omega) \;=\; 
\sum_{k=1}^P a_k(\omega)\,f_k(s),
\end{equation}
where the coefficients $a_k(\omega)$ depend on type of response surface. The approximation can be cheaply sampled by computing the coefficients $a_k(\omega)$ at the $M$ points in $\Omega$ needed to estimate $q_s$. This approximation introduces additional error in the estimate of $q_s$, and one should validate that $\tilde{f}$ is sufficiently accurate---e.g., that $P$ is large enough to produce a good approximation. Asymptotically, an $L_2$-convergent response surface implies convergence in distribution, i.e., the pdf of $\tilde{f}$ converges to the pdf of $f$~\cite[Chapter 2]{durrett2005}. However, this well-known result does not account for the finite sampling used to estimate the pdfs, and asymptotic results do not always give confidence in the case when $P$ is small. 

The airfoil example in section \ref{sec:airfoil} has a smooth response that is a function of one parameter representing uncertainty, and the system uses a relatively expensive CFD solver in two spatial dimensions. We use a response surface---validated several points in the design space---for both the response and its partial derivatives with respect to the design variables. 

\subsection{Implementation details}
\label{sec:imp}

\noindent In this section, we collect the pieces needed to implement a numerical solver for the optimization \eqref{eq:opt}. This summarizes the method and provides some details about the choices we make in implementation.

\paragraph{Optimization package} From MATLAB's Optimization Toolbox, we use the \texttt{fmincon} function with \texttt{Algorithm} option set to \texttt{sqp} (sequential quadratic program). An open source alternative is  SciPy's \texttt{minimize} function from its Optimization package with the \texttt{method} set to \texttt{SLSQP} (sequential least-squares quadratic program). We provide subroutines for computing the objective function, implemented as \eqref{eq:aopt} with the approximation \eqref{eq:vqs}, and the objective's gradient, implemented as \eqref{eq:equation_gradients}. 

\paragraph{Numerical integration} We hand code a trapezoidal rule to evaluate the objective function \eqref{eq:aopt}. This includes forming the the diagonal matrix $\mW$ in \eqref{eq:aopt} and \eqref{eq:equation_gradients}. 

\paragraph{Kernel density estimation} From MATLAB's Statistics and Machine Learning Toolbox, we use the \texttt{ksdensity} function for kernel density estimation. The default kernel is the Gaussian as in \eqref{eq:gaussk}, and the interface takes an optional bandwidth parameter argument, which we use to implement the heuristic described in section \ref{sec:kde}. We hand code the gradient of the Gaussian kernel to compute $\mK'$ in \eqref{eq:equation_gradients}. An open source alternative is SciPy's \texttt{gaussian\_kde} function, which also accepts a user-specified bandwidth parameter. 

\paragraph{Response surfaces} The variety of response surface types and application scenarios make it difficult to create a general purpose toolbox for response surfaces. For the specific experiment in this paper, we use a hand-coded least-squares-fit 5th degree global polynomial approximation for the univariate response. Since the airfoil application has only a single parameter representing uncertainty, we are able to visually validate the response surface quality with two-dimensional plots.

%% file: examples.tex
\section{Numerical examples}
\label{sec:examples}

\noindent We study two numerical examples to explore the characteristics of the proposed density-matching OUU approach. The first is a simple response function that produces surprisingly complex behavior in the optimization. The second is an airfoil shape optimization problem with uncertainty in the freestream mach number. 

\subsection{Simple response function}
\label{sec:simple}

\begin{figure}
\centering
\includegraphics[natwidth=7.78in, natheight=5.83in, width=15cm]{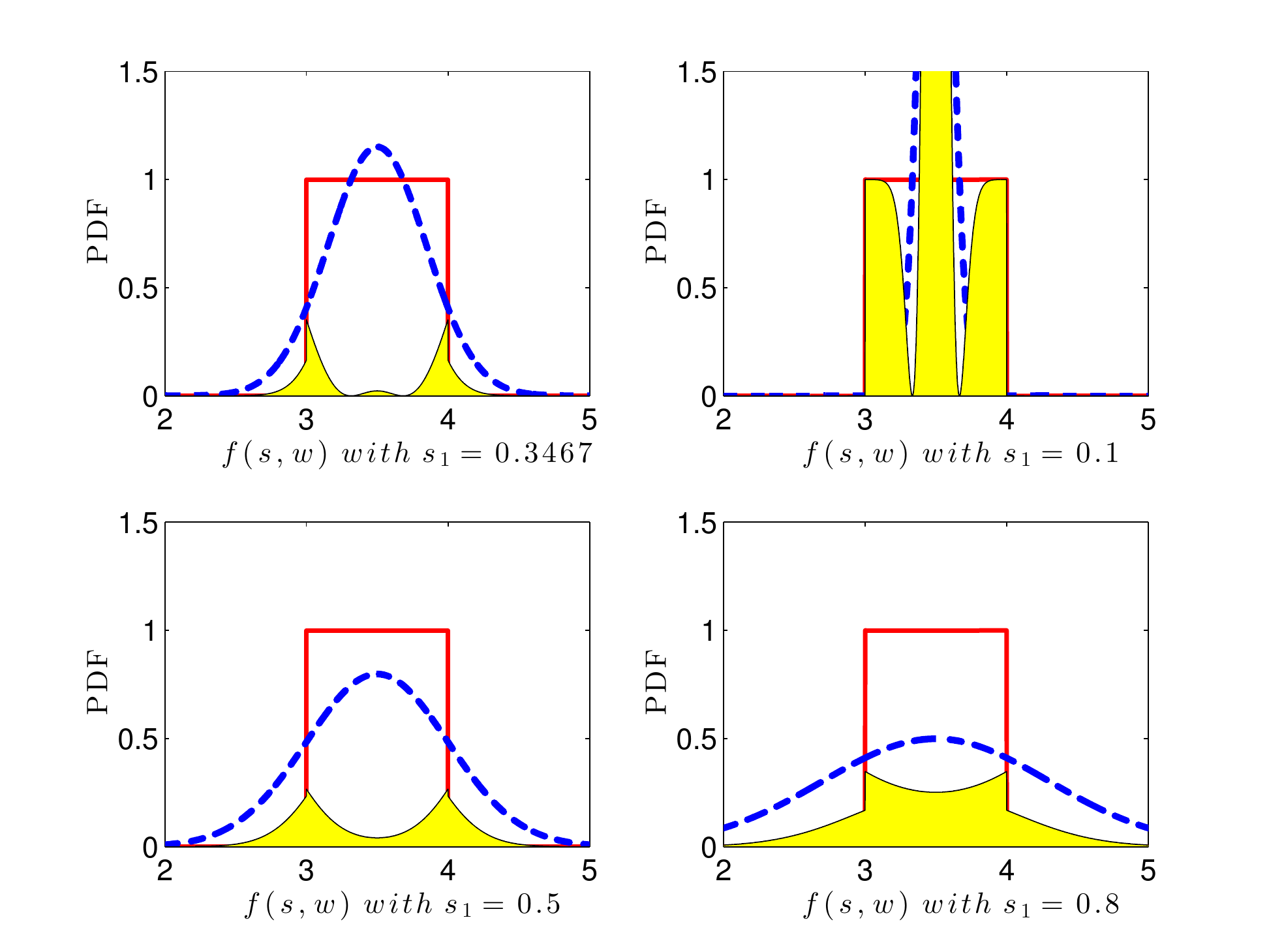}
\caption{Response \eqref{eq:simple} with different values of the design parameter $s$. The red rectangle is the uniform target pdf and the blue dashed line is the response pdf. The goal is to minimize the yellow area in the figures.}
\label{linearmodel}
\end{figure}

\noindent Consider the model
\begin{equation}
\label{eq:simple}
f(s,\omega) \;=\; s\omega+3.5,
\end{equation}
where $s\in\mathbb{R}$ and $\omega$ is a standard normal random variable, so $f$ is a normal random variable with mean 3.5 and standard deviation $s$. The goal is to find $s$ that minimizes the $L_2$-norm distance between the pdf of $f(s,w)$ and a uniform target pdf 
\begin{equation}
t(f) \;=\; \left\{
\begin{array}{ll}
1, & \mbox{ for $f\in[3,4]$,}\\
0, & \mbox{ elsewhere.} 
\end{array}
\right.
\end{equation}
The pdf of $f$ is
\begin{equation}
q_s(f) \;=\; \frac{1}{\sqrt{2\pi} s}\exp\left(\frac{-\left(f-3.5\right)^{2}}{2s^2}\right).
\end{equation}
The objective function is
\begin{equation}
\label{equ1}
\begin{aligned}
d(t,q_{s}) &=\int_{-\infty}^{\infty}\left(t(f)-q_s(f)\right)^{2}\,df\\
&=\int_{-\infty}^{3} q_s(f)^{2} \,df 
+ \int_{4}^{\infty} q_s(f)^{2}\,df 
+ \int_{3}^{4}\left(1-q_s(f)\right)^{2}\,df.
\end{aligned}
\end{equation}
Figure \ref{linearmodel} plots the two densities and their squared difference. The red density is the target pdf and the blue dashed line is the response pdf. The yellow-shaded region is the squared difference in \eqref{equ1}. Figure \ref{linearmodel} shows the effect of varying $s$ on the yellow region, with $s=0.3467$ yielding the smallest area under the curve, minimizing the distance \eqref{equ1}. Note the complexity of the integrand. Table~\ref{linearmodel1} shows the effect of varying the number of quadrature points. Above 1000 points the error in the optimal $s$ is in the fourth decimal place. Table~\ref{linearmodel2} repeats this study with kernel density estimates for the response pdf. 

\begin{table}
 \begin{center}
   \caption{Effect of the number of quadrature points for the simple problem \eqref{eq:simple}.}
   \label{linearmodel1}
   \begin{tabular}{lcc}
%%   \toprule
    Number of quadrature points & Optimal $s$ & Minimum distance \\
    \hline
   10 & 0.3531 & 0.00670\\  
   100 & 0.3556 & 0.017385 \\  
   1000 & 0.3469 & 0.016048 \\
   10000 & 0.3467 & 0.016022 \\
   100000 & 0.3467 & 0.016023 \\
%%   \toprule
   \end{tabular}
 \end{center}
\end{table}

\begin{table}
 \begin{center}
   \caption{Effect of the number of quadrature points and kernel density estimate (KDE) samples for the simple problem \eqref{eq:simple}.}
   \label{linearmodel2}
   \begin{tabular}{lccc}
%%   \toprule
    Number of quadrature points & Optimal $s$ & Minimum distance & Number of samples \\
    \hline
   10 & 0.3488 & 0.006759 & $10^5$\\  
   10 & 0.3505 & 0.006729 & $10^6$\\  
   100 & 0.3513 & 0.01743 & $10^5$ \\ 
   100 & 0.3537 & 0.01741 & $10^6$ \\
   1000 &  0.3464 & 0.01610 & $10^5$ \\
   1000 &  0.3465 & 0.01610 & $10^6$ \\
   10000 & 0.3466 & 0.01608 & $10^5$ \\
   10000 & 0.3467 & 0.01605 & $10^6$ \\
%%   \toprule
   \end{tabular}
 \end{center}
\end{table}

\subsection{Airfoil design}
\label{sec:airfoil}
\noindent Next we apply the density matching scheme to the design of an airfoil under uncertainty. MATLAB and Python codes used for this numerical study can be found at \url{https://github.com/psesh/density-matching}. The airfoil used in this example is a NACA0012 at a Reynolds number of $10^6$ and an angle of attack of $5^{\circ}$. The uncertainty is in the inlet Mach number, which is characterized by a $\beta(2,2)$ distribution between Mach numbers of 0.66 and 0.69. Flow computations for this airfoil are carried out by solving the compressible Euler equations using Stanford University's SU2 flow solver~\cite{SU2}. The airfoil is parameterized with 16 Hicks-Henne bump functions: 8 on the upper surface and 8 on the lower surface. The design space is the height of each bump; a point in the design space produces a perturbation from the NACA0012 shape. The height ranges and locations for the bumps are shown in Table~\ref{hicks}. 

\begin{table}
 \begin{center}
   \caption{Hicks-Henne bump function heights and locations as a proportion of chord. The heights and locations are shown for the upper surface; the lower surface has the same parameterization.}
   \label{hicks}
   \begin{tabular}{lcc}
%%   \toprule
    Location & Bump amplitude \\
    \hline
0.05 & $\pm 0.0007$\\  
0.15 & $\pm 0.0030$\\  
0.30 & $\pm 0.0090$\\  
0.40 & $\pm 0.0090$\\  
0.55 & $\pm 0.0090$\\  
0.65 & $\pm 0.0060$\\     
0.75 & $\pm 0.0030$\\  
0.90 & $\pm 0.0007$\\ 
%%   \toprule
   \end{tabular}
 \end{center}
\end{table}

For a point in the design space, the airfoil mesh is deformed using a torsional spring analogy. The flow solver runs on the new mesh producing the lift-to-drag ratio $L/D$, which is the response of interest. To connect to the notation in Sections \ref{sec:math} and \ref{sec:discretization}, the response $f$ is $L/D$, the design variables $s$ are the 16 Hicks-Henne bump heights, and the random variable $\omega$ is the Mach number with a $\beta(2,2)$ density on the interval $[0.66,0.69]$.

\subsection{Robust design optimization}
\label{sec:rdo}
\noindent We begin our investigation with a multi-objective RDO. There are three aims of this exercise:
\begin{enumerate}
\item estimate the cost of an RDO for an unconstrained 16 parameter problem,
\item illustrate some limitations of RDO objectives,
\item generate a Pareto front to qualitatively compare with density-matching results.
\end{enumerate}
The RDO problem is
\begin{equation}
\label{rdo2}
\underset{s\in\sS}{\operatorname{minimize}} \;\; (\mathbb{E}\left[ L/D \right])^{-1} \mbox{ and } \operatorname{Var}\left[L/D\right],
\end{equation}
where $\mathbb{E}[\cdot]$ is the mean and $\operatorname{Var}[\cdot]$ is the variance. These two moments both depend on the design variables $s$ that parameterize the airfoil shape. We use the genetic algorithm NSGA-II~\cite{NSGA} to estimate the Pareto front for \eqref{rdo2}. The default parameters for NSGA-II are shown in Table~\ref{nsga}. 
\begin{table}
 \begin{center}
   \caption{NSGA-II algorithm parameters used for RDO}
   \label{nsga}
   \begin{tabular}{lll}
%%   \toprule
    Parameters & Value \\
    \hline
Population size & 100 \\
Number of generations & 35 \\
Crossover probability & 0.9 \\
Mutation probability & 0.0625 \\
Crossover distribution index & 20\\
Mutation distribution index & 20\\
%%   \toprule
   \end{tabular}
 \end{center}
\end{table}

We use a population size of 100 with 35 generations, yielding a total of 3500 function calls from the optimizer. Each function call uses 21 CFD computations to fit a least-squares 5th degree polynomial response surface, which is used to estimate the objectives in \eqref{rdo2}. Thus, the RDO study used a total of $3500 \times 21=73,500$ CFD runs.

\begin{figure}
\centering
\includegraphics[natwidth=792, natheight=612, width=10.5cm]{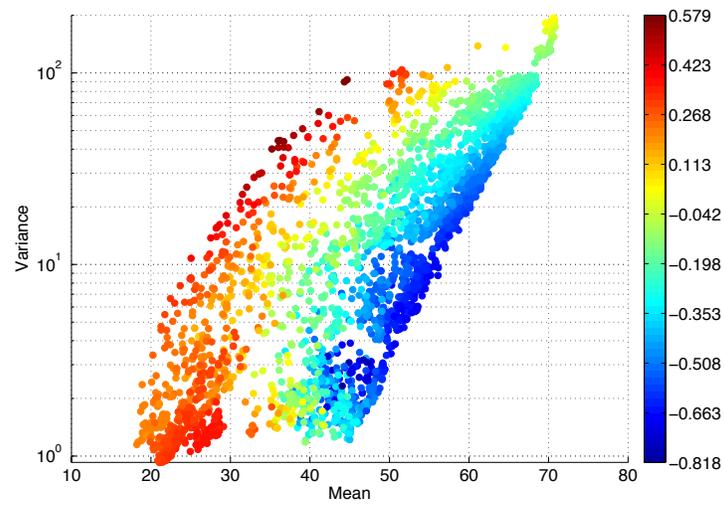}
\caption{Results the RDO airfoil design problem. Mean and log-variance of $L/D$ are plotted on the horizontal and vertical axes, respectively. Individual designs are colored by their skewness values as shown on the color bar legend.}
\label{front}
\end{figure}

Figure \ref{front} plots the moments from each computation. The mean is plotted on the horizontal and the variance on the vertical on a logarithmic scale. The nominal NACA0012 design has a mean $L/D$ ratio of 27.2356 and a variance of 7.4615. Figure \ref{front} indicates the skewness values of individual designs by the marker color. Negatively skewed designs are blue while positively skewed designs are red. The RDO took approximately two days to run on an 8-core workstation.

\subsection{Density-matching with a designer-specified target}
\label{sec:adnaca}
\begin{figure}
\centering
\includegraphics[natwidth=7.78in, natheight=5.83in, width=8cm]{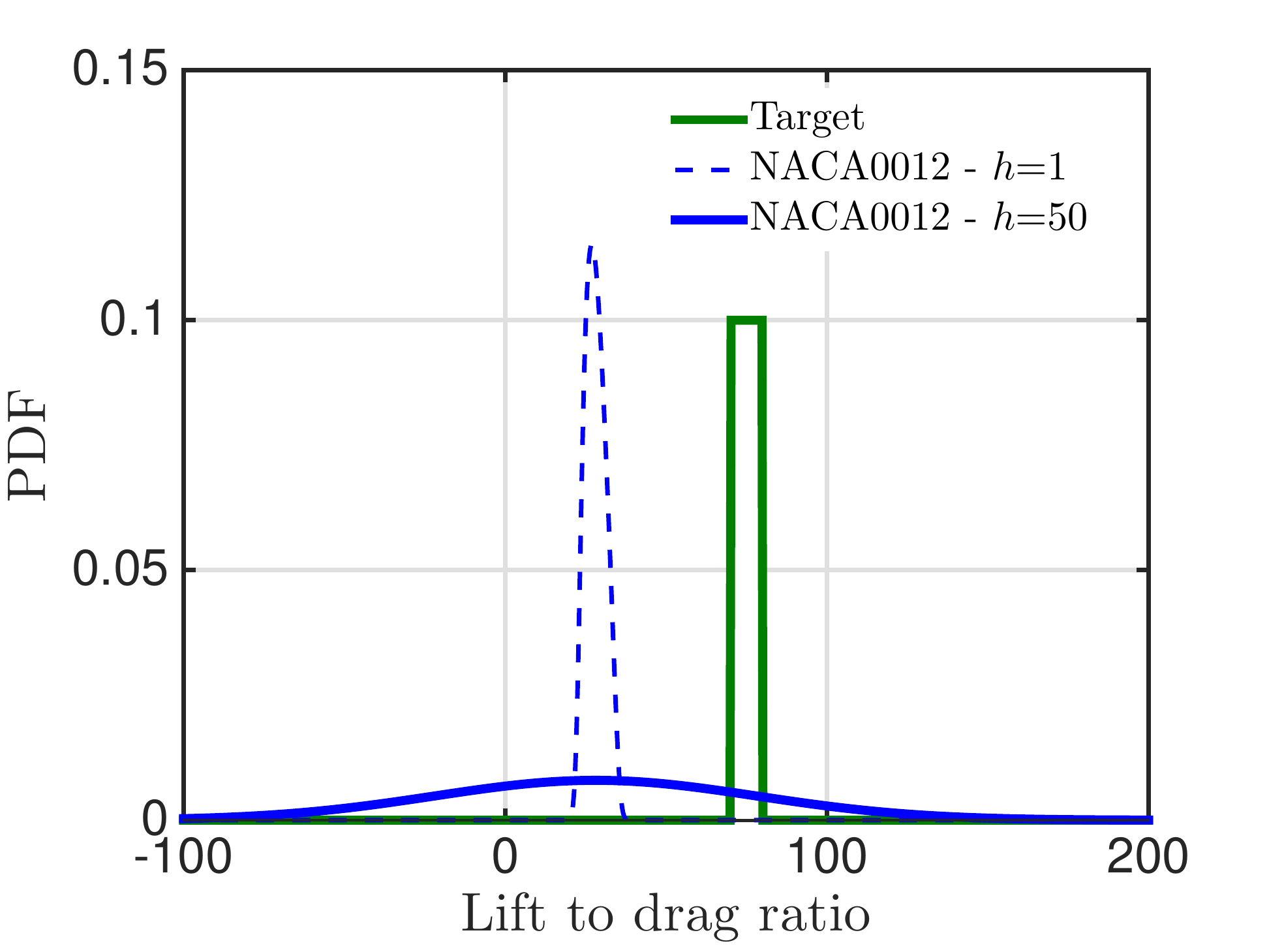}
\caption{Initial (NACA0012) design pdfs and the uniform target. The initial design is shown with bandwidth parameter values of 1 and 50.}
\label{overlap}
\end{figure}

\noindent We begin with some details of the density-matching optimization applied to the airfoil design. The initial design for all cases is the NACA0012 airfoil. To compute the response pdf's kernel density estimates, we draw $10^5$ independent samples from a least-squares-fit 5th degree polynomial response surface of $L/D$ as a function of Mach number. We fit the response surface with $P=21$ flow computations at uniformly spaced Mach numbers between 0.66 and 0.69; the 21 independent runs were executed in parallel. The nominal design (the NACA0012 airfoil) produces an estimated response pdf such that the three targets reside in the tails; see Figure \ref{overlap}. We use the two-stage heuristic discussed in Section \ref{sec:kde}. We ran 3 iterations of the SQP solver with a large bandwidth of $h=50$. The remaining iterations used a bandwidth $h=1$---slightly larger than $h_{\text{opt}}$ in \eqref{eq:scott}.
We use a trapezoidal rule with $N=2500$ quadrature points on the interval $[-100,150]$ in the objective function \eqref{eq:aopt}. The negative lower bound on $L/D$ accommodates kernel density estimates with large bandwidths. We tested larger values of $N$, but they did not lead to substantial improvement in the optimized design.

\begin{figure}
\centering
\includegraphics[natwidth=962, natheight=561, width=13.5cm]{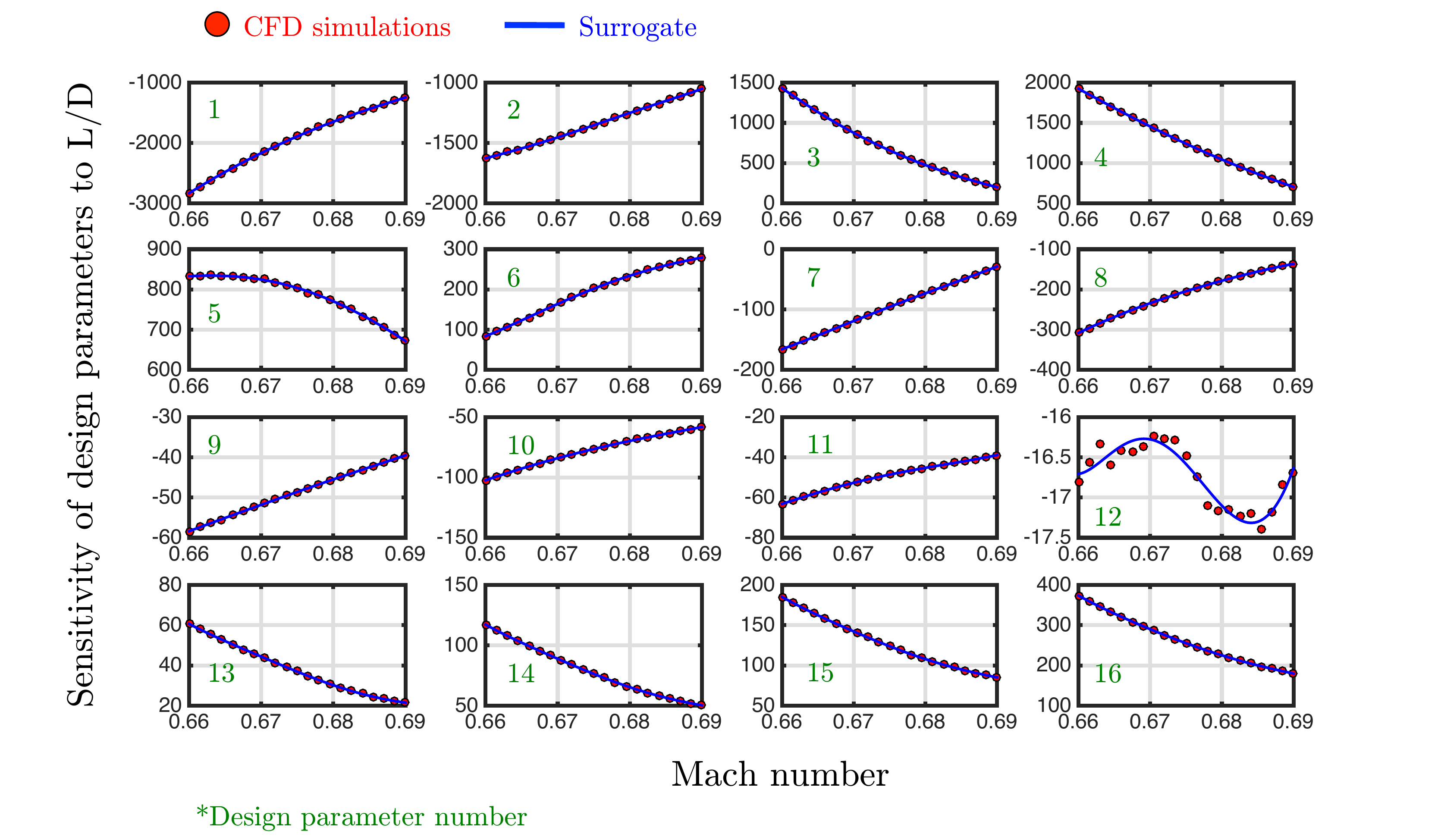}
\caption{Polynomial response surfaces and adjoint-based partial derivatives of $L/D$ with respect to Mach number for a random design point.}
\label{adjsens}
\end{figure}

The gradient in \eqref{eq:equation_gradients} includes the matrix of kernel evaluations $\mK$, their derivatives $\mK'$, and design parameter sensitivities in the matrix $\mF'$. Entries of $\mK$ and $\mK'$ are computed using the Gaussian kernel \eqref{eq:gaussk}. The size of these matrices are determined by the number of quadrature points and the number of random samples used to estimate the pdfs. Both $\mK$ and $\mK'$ have dimensions $2500 \times 10^5$. We used SU2's adjoint capabilities at the $P=21$ Mach numbers to compute partial derivatives of $L/D$ with respect to 16 bump heights. We used a least-squares-fit 5th degree polynomial response surface to approximate the partial derivatives at all the sample points needed to compute the elements of $\mF'$ in \eqref{eq:grads}. The response surfaces for all 16 partial derivatives are shown in Figure \ref{adjsens} for a perturbed NACA0012. From these and other similar plots, we determined that the polynomial response was sufficiently accurate. We used least-squares to avoid interpolating noisy partials, such as parameter 12 in Figure \ref{adjsens}. Table~\ref{aggressive_study} summarizes the parameters used in the optimization. 
\begin{table}[h]
 \begin{center}
   \caption{Density-matching parameters}
   \label{aggressive_study}
   \begin{tabular}{lll}
%%   \toprule
    Quantity & Definition & Value \\
    \hline
$N$  & quadrature points & 2500 \\
$M$ & random samples & $10^5$ \\
$n$ & design parameters &  16 \\
$h_{\text{stage 1}}$ & stage 1 bandwidth parameter & 50.0\\
$h_{\text{stage 2}}$ & stage 2 bandwidth parameter & 1.0\\
$f_\ell$ & lower bound for $L/D$ & -100 \\
$f_u$ & upper bound for $L/D$ & 150 \\
$K_h$ & kernel function & Gaussian
%%   \toprule
   \end{tabular}
 \end{center}
\end{table}

We assume that the designer has provided us with three target pdfs for the response $L/D$:
\begin{enumerate}
\item a uniform density in the interval $[75,80]$,
\item a Gaussian density with mean $50$ and variance $10$,
\item a $\beta(1.5,3.5)$ density on the interval $[50,80]$.
\end{enumerate}
In what follows, we describe the results for each of these three cases. 

\subsubsection{Uniform target}
\begin{figure}
\begin{subfigmatrix}{2}% number of columns
\subfigure[]{\includegraphics[natwidth=7.78in, natheight=5.83in, width=5.5cm]{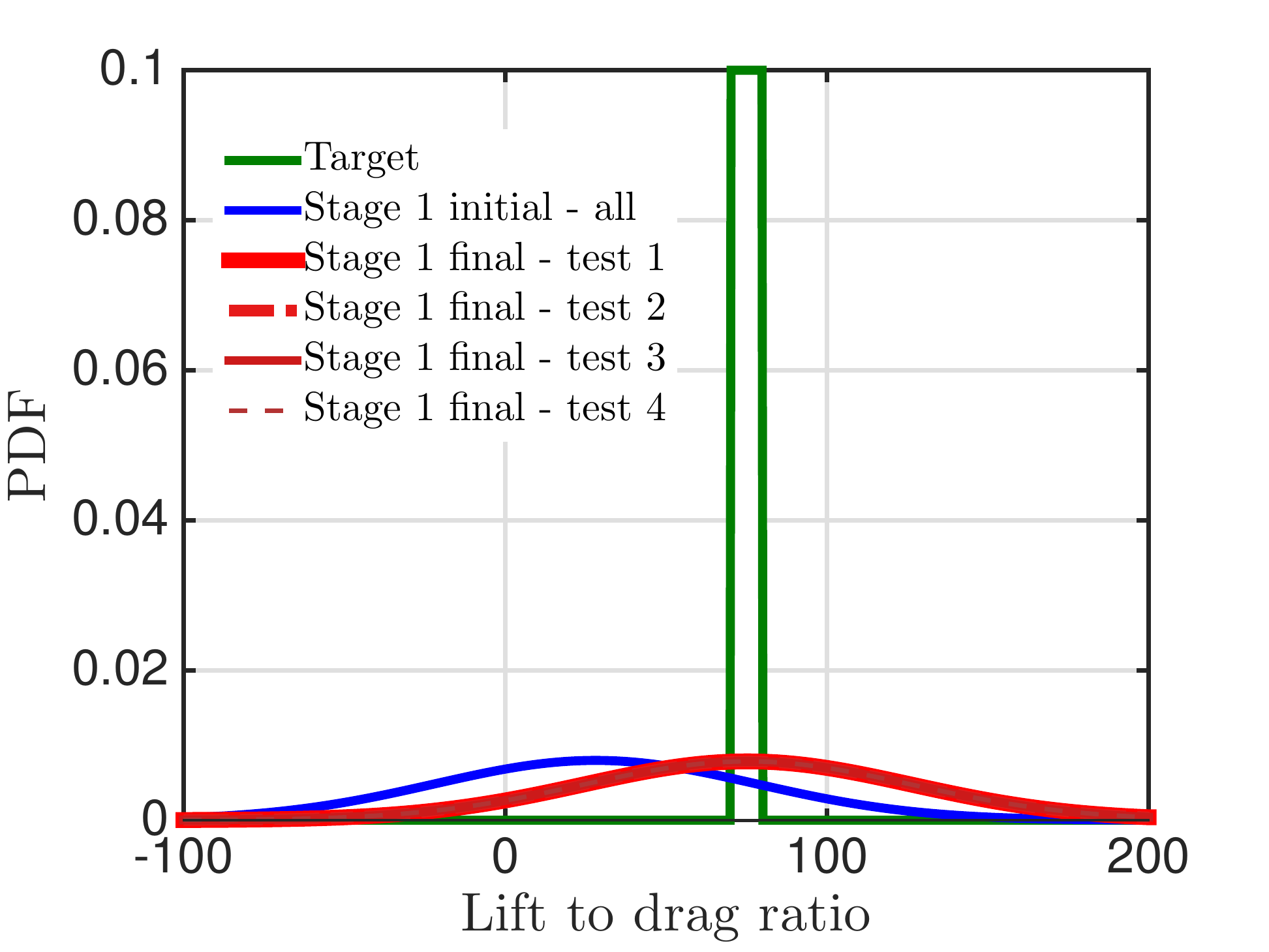}\label{fig:a1}}
\subfigure[]{\includegraphics[natwidth=7.78in, natheight=5.83in, width=5.5cm]{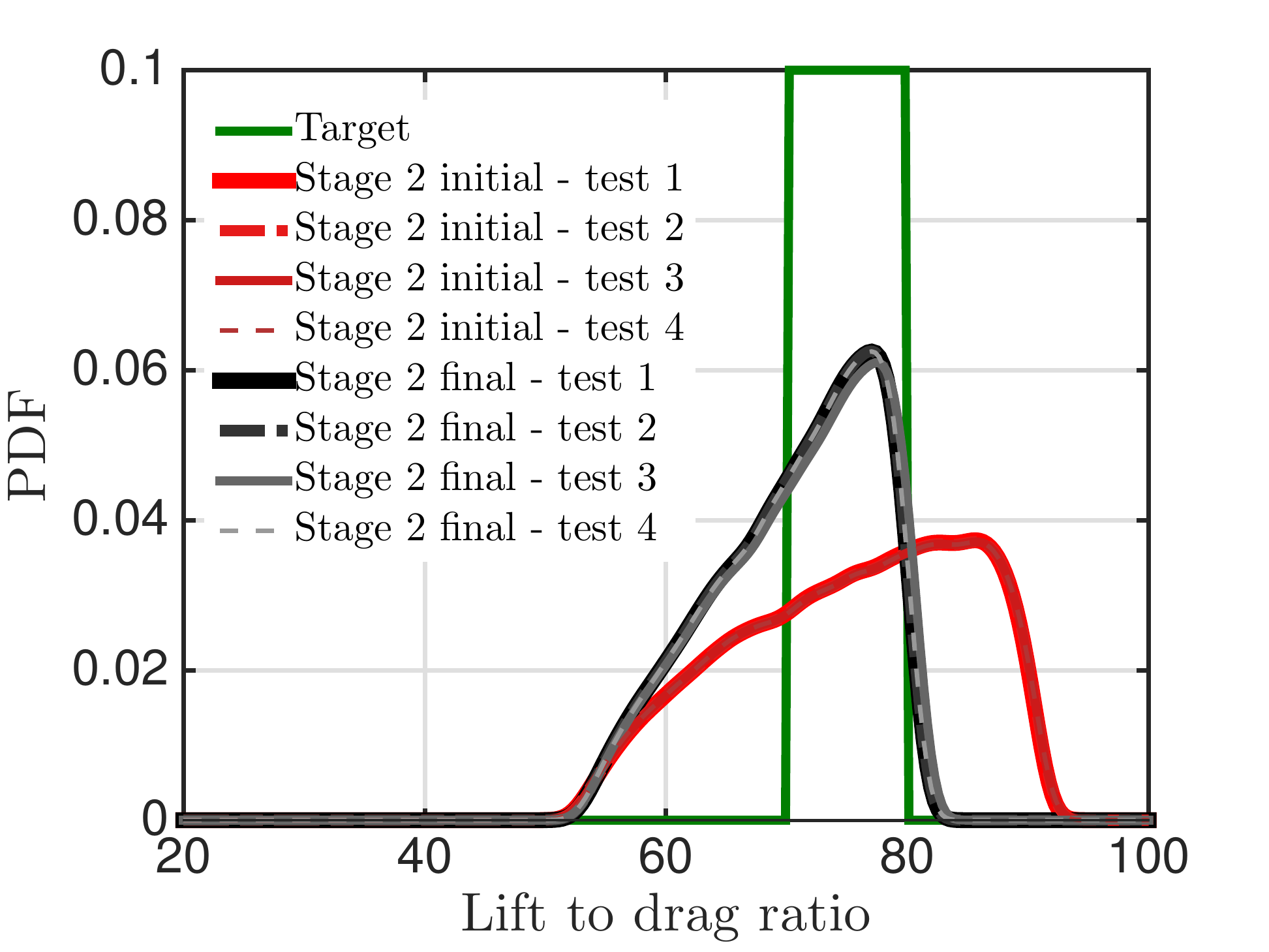}\label{fig:a2}}
\subfigure[]{\includegraphics[natwidth=7.78in, natheight=5.83in, width=5.5cm]{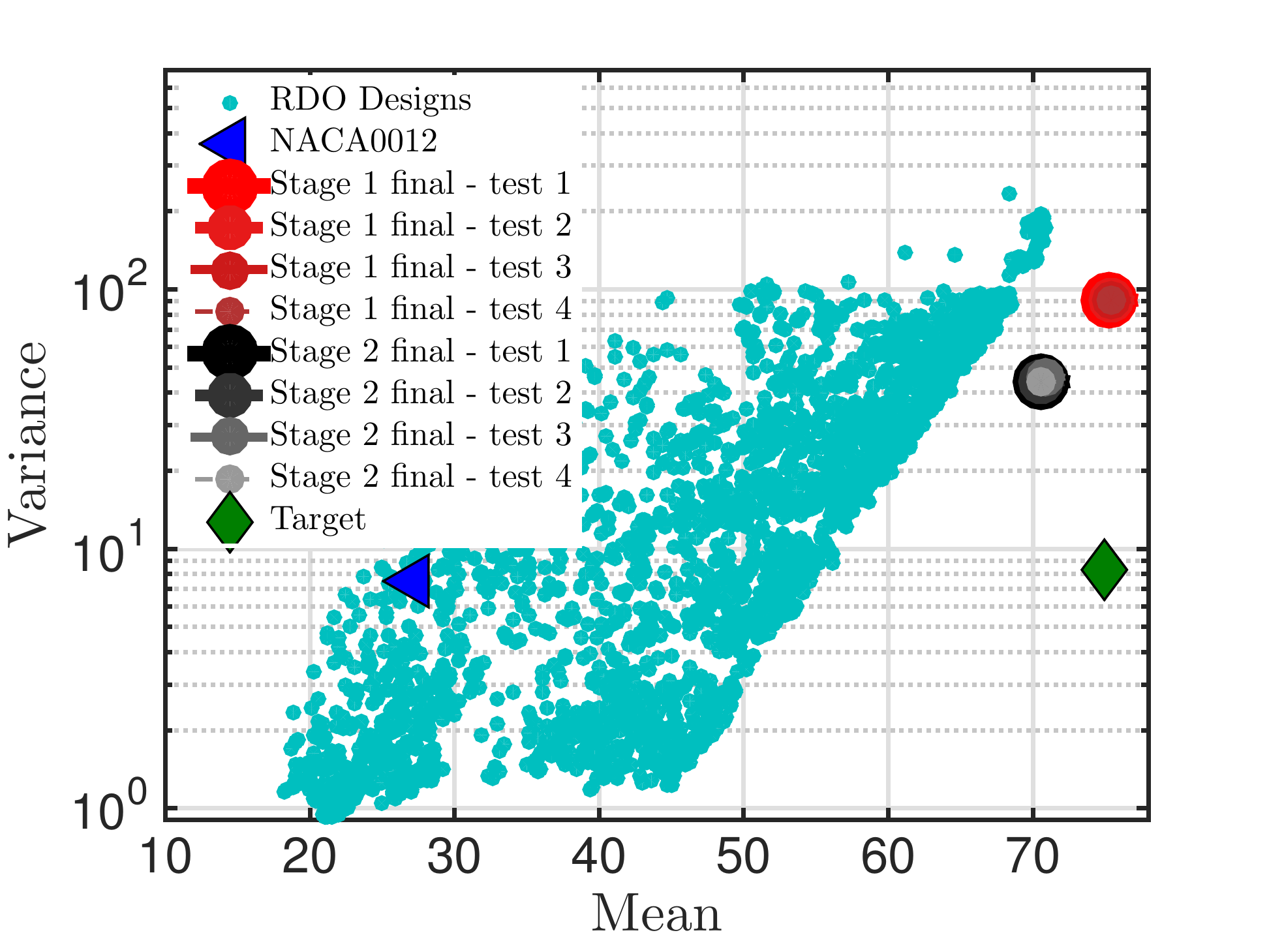}\label{fig:a3}}
\subfigure[]{\includegraphics[natwidth=7.78in, natheight=5.83in, width=5.5cm]{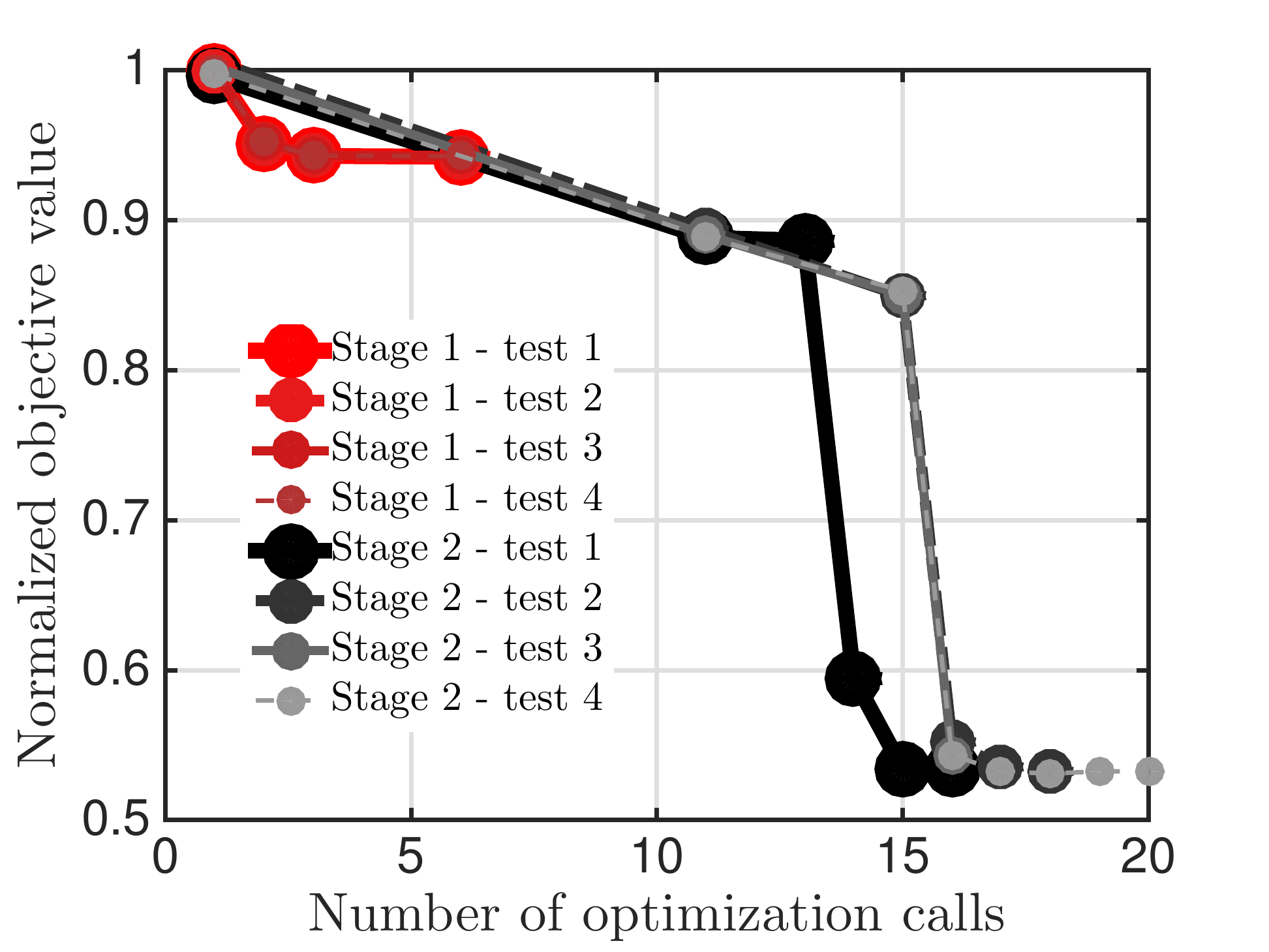}\label{fig:a4}}
\subfigure[]{\includegraphics[natwidth=7.78in, natheight=5.83in, width=5.5cm]{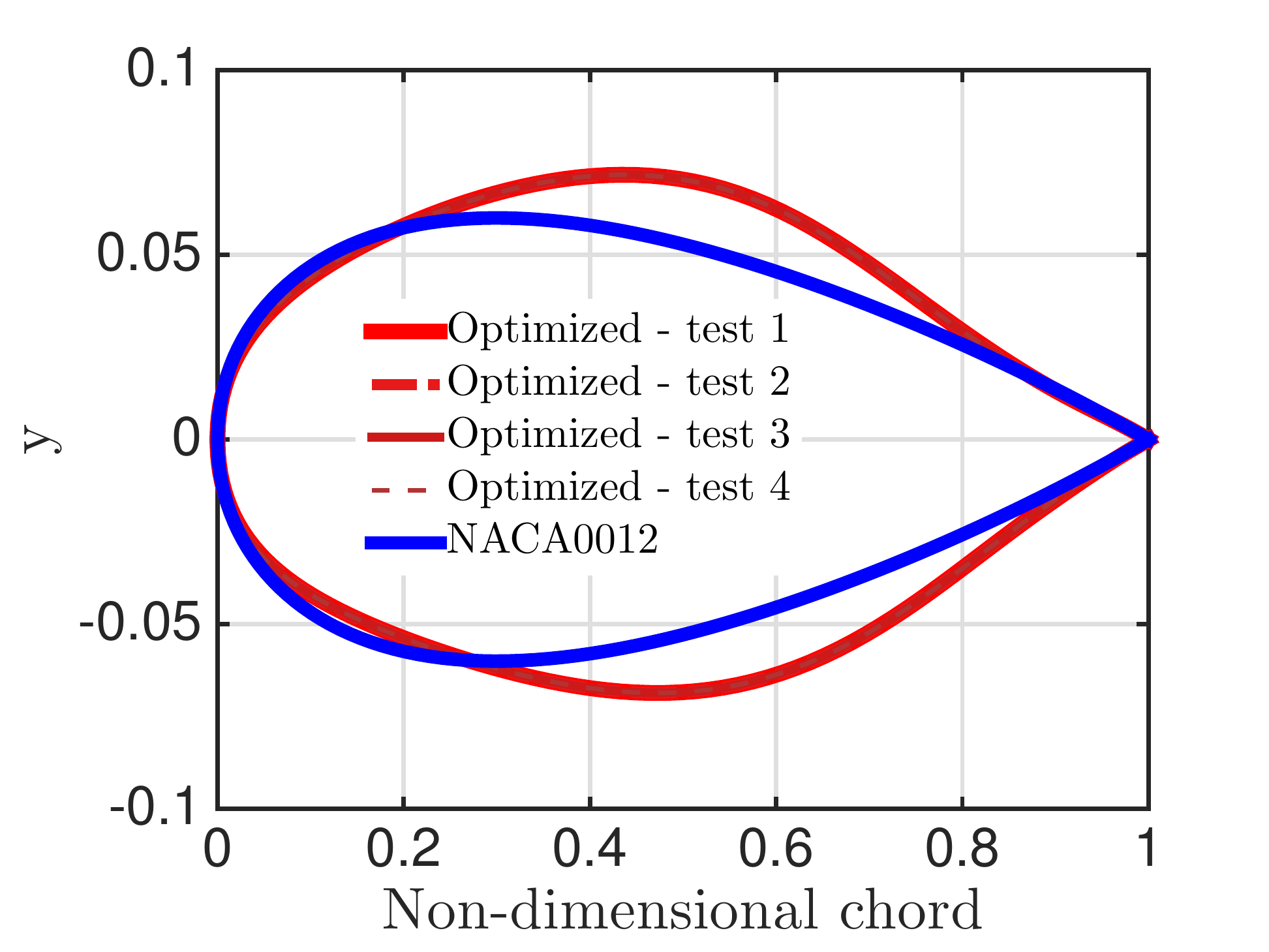}\label{fig:a5}}
\subfigure[]{\includegraphics[natwidth=7.78in, natheight=5.83in, width=5.5cm]{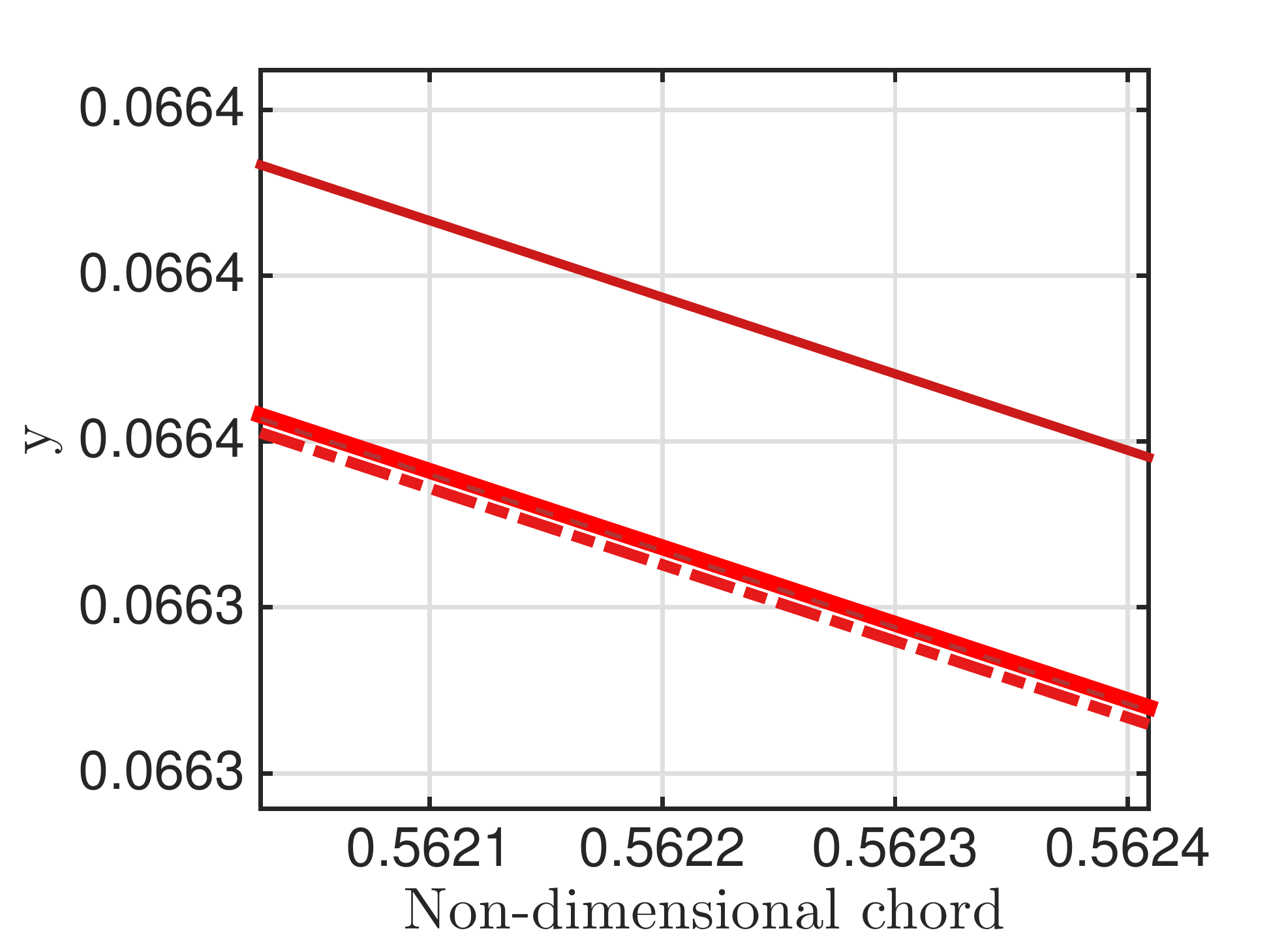}\label{fig:a6}}
\end{subfigmatrix}
\caption{Uniform target results: (a) stage 1 ($h=50$), (b) stage 2 ($h=1$), (c) comparison with RDO designs, (d) stage 1 and stage 2 convergence plots, (e) stage 2 optimal designs with a close-up in (f).}
\label{A}
\end{figure}

\noindent For the uniform target, we repeat the density-matching optimization four times to see the effects of the random sampling used to estimate the response pdf at each optimization iteration. The four independent trials produced nearly identical results. The results with the uniform target are shown in Figure \ref{A}. Figure \ref{fig:a1} compares the kernel density estimate for the initial design with bandwidth $h=50$ (blue), the uniform target (green), and the final designs from stage 1 of the optimization using bandwidth $h=50$ across all four independent trials (red). Using the final stage 1 designs with bandwidth $h=1.0$ produces the red density estimates shown in Figure \ref{fig:a2}; these are the initial designs for stage 2 using bandwidth $h=1$. The final optimized designs from stage 2 produce the black density estimates in Figure \ref{fig:a2}. We plot the final designs from both stages in Figure \ref{fig:a3} to compare them with the RDO designs obtained in section \ref{sec:rdo}. In this case, the results from both stage 1 and stage 2 lie well beyond the Pareto front obtained from NSGA-2. 

Optimization convergence histories for all four tests are shown in Figure \ref{fig:a4}. All objective function values are normalized by the initial objective value. Stage 1 for this target took an average of 1 hour and 34 minutes while stage 2 took an average of 3 hours and 36 minutes. While there are some differences in the convergence histories, the differences in the final designs are very minor across the four independent trials. The final designs from stage 2 are plotted in Figure \ref{fig:a5} with a close-up in Figure \ref{fig:a6}. 

\subsubsection{Gaussian target}
\begin{figure}
\begin{subfigmatrix}{2}% number of columns
\subfigure[]{\includegraphics[natwidth=7.78in, natheight=5.83in, width=5.5cm]{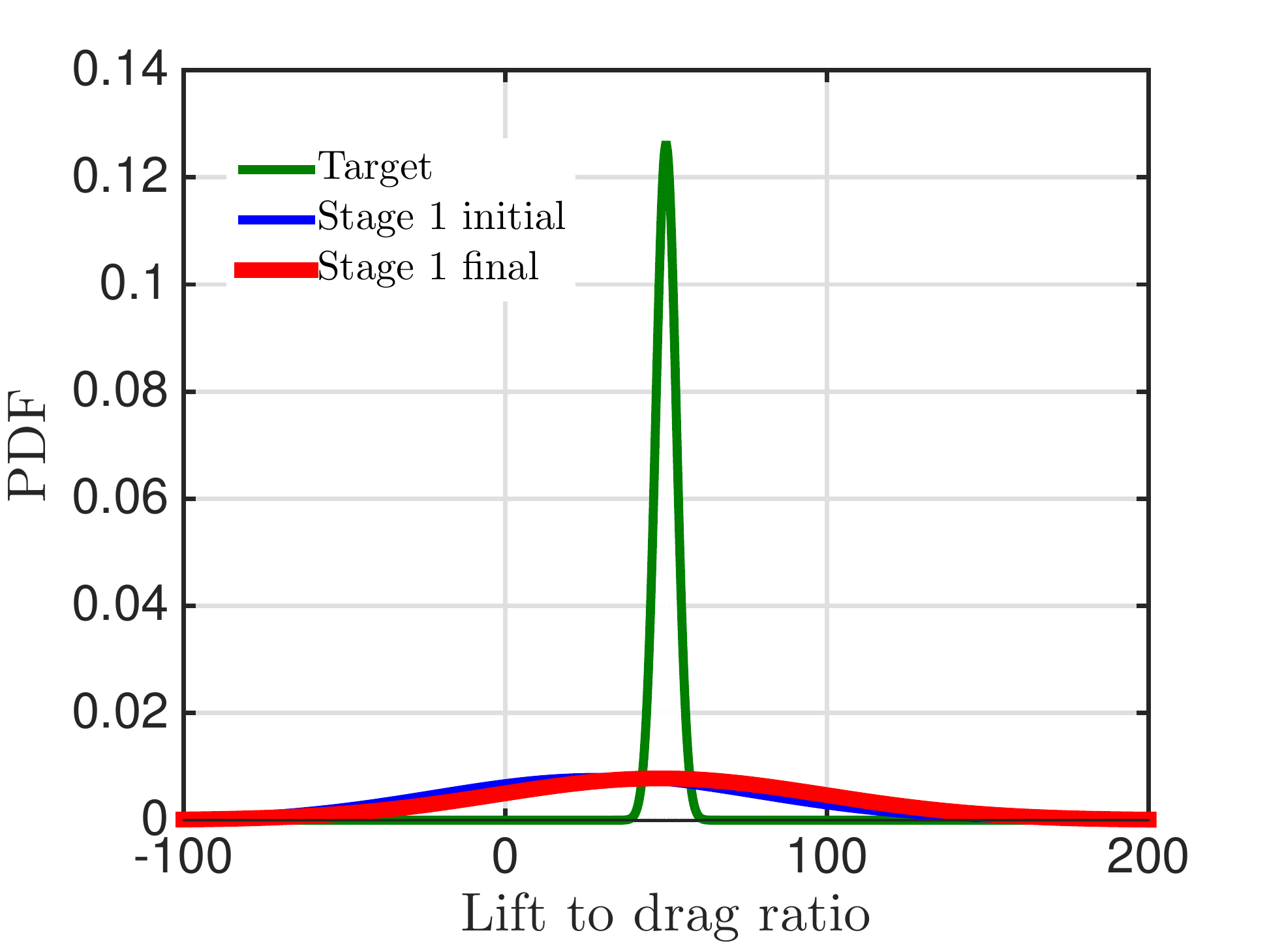}\label{fig:b1}}
\subfigure[]{\includegraphics[natwidth=7.78in, natheight=5.83in, width=5.5cm]{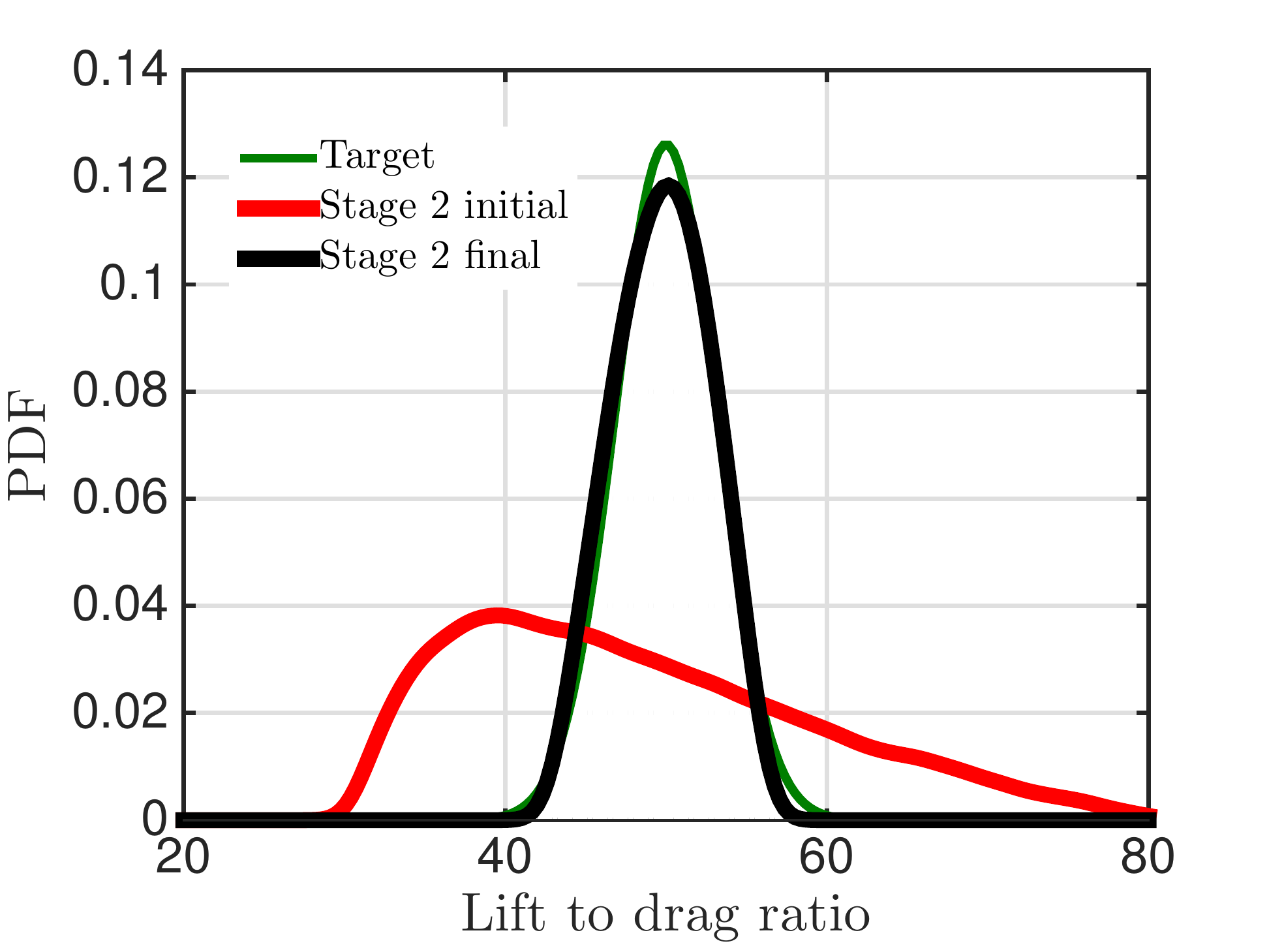}\label{fig:b2}}
\subfigure[]{\includegraphics[natwidth=7.78in, natheight=5.83in, width=5.5cm]{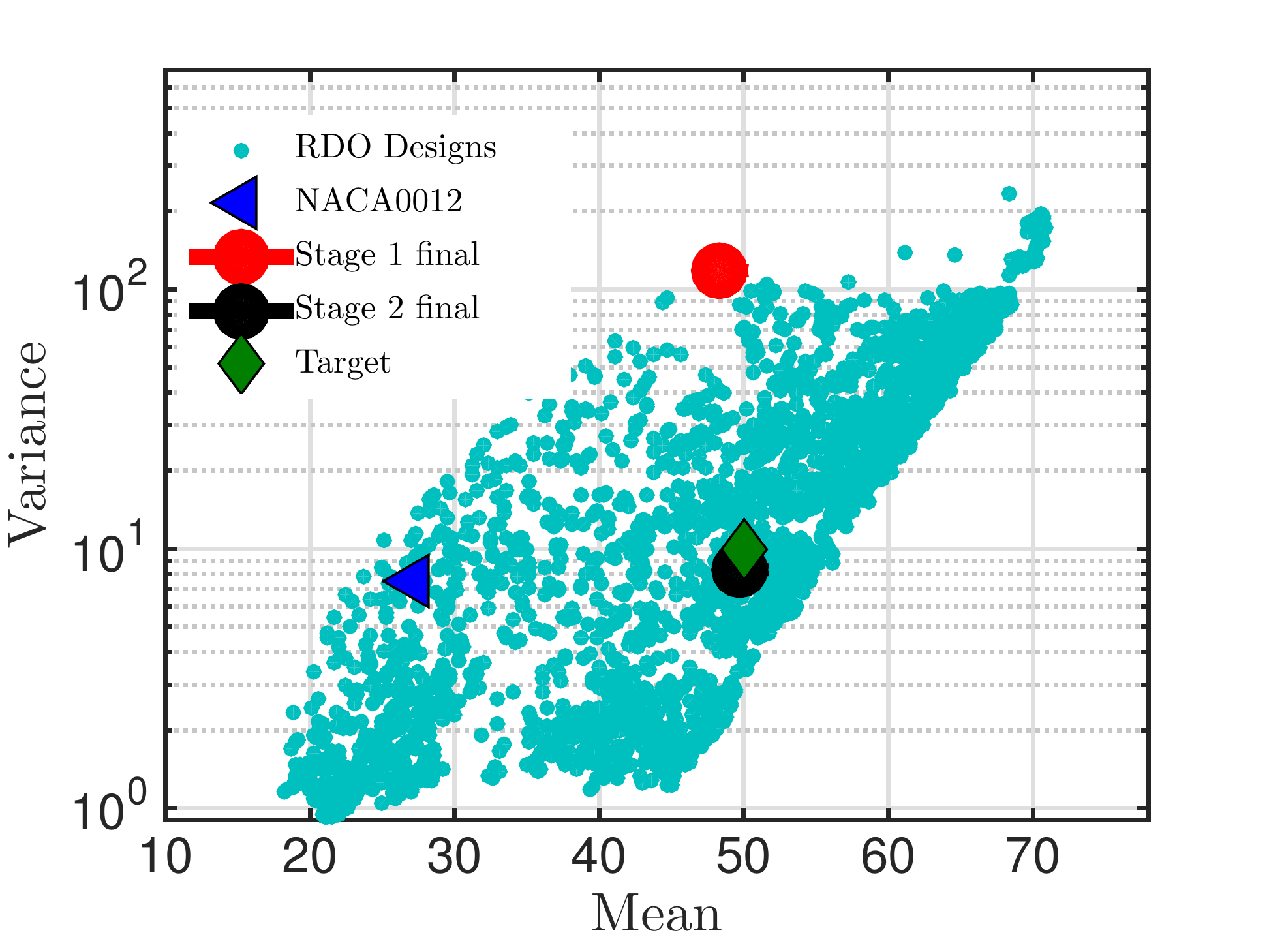}\label{fig:b3}}
\subfigure[]{\includegraphics[natwidth=7.78in, natheight=5.83in, width=5.5cm]{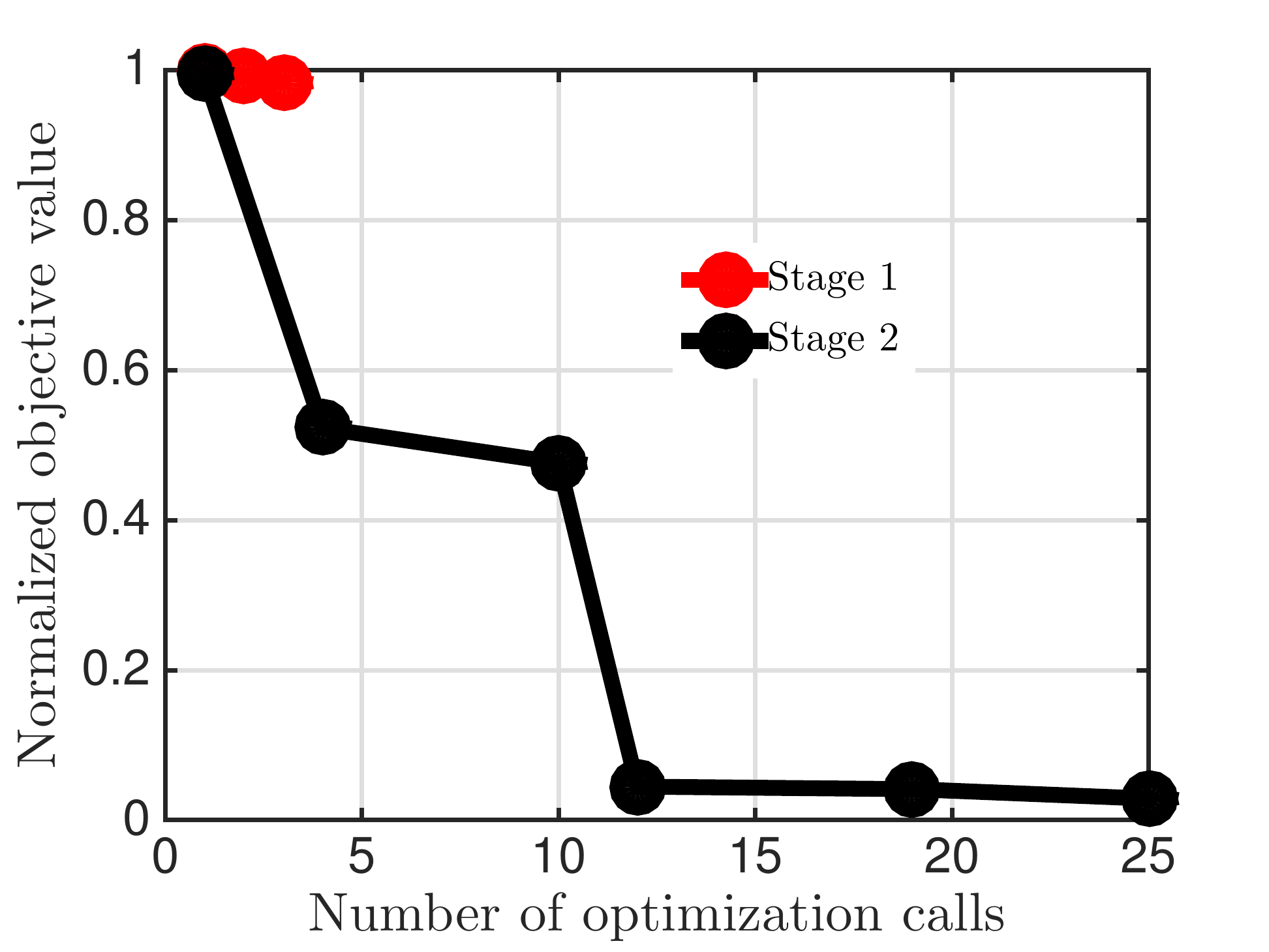}\label{fig:b4}}
\subfigure[]{\includegraphics[natwidth=12.61in, natheight=5.83in, width=8cm]{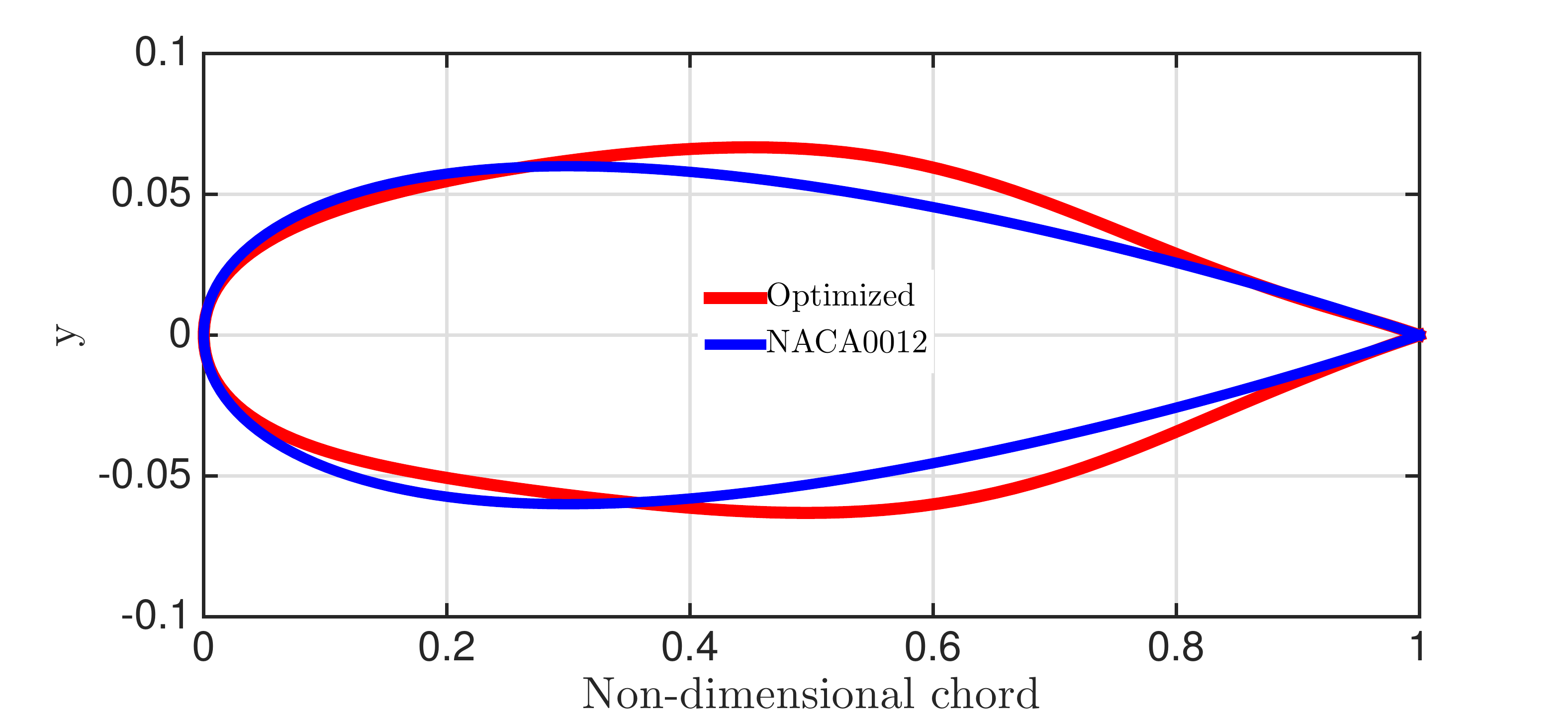}\label{fig:b5}}
\end{subfigmatrix}
\caption{Gaussian target results: (a) stage 1 ($h=50$), (b) stage 2 ($h=1$), (c) comparison with RDO designs, (d) stage 1 and stage 2 convergence plots, (e) stage 2 optimal design.}
\label{Gauss}
\end{figure}

\noindent Figure \ref{Gauss} shows results with the Gaussian target similar to Figure \ref{A}. In this case, we ran a single test instead of four independent trials. Figure \ref{fig:b1} shows the initial and final designs from stage 1, and Figure \ref{fig:b2} shows the initial and final designs for stage 2. Note that the final design for stage 1 is the initial design for stage 2. The final design is extremely close to matching the target. Figure \ref{fig:b4} shows that changes in the objective function value are negligible after 12 optimization calls in the second stage. Here stage 1 took 46 minutes while stage 2 took 5 hours and 12 minutes.

\subsubsection{Beta target}
\begin{figure}
\begin{subfigmatrix}{2}% number of columns
\subfigure[]{\includegraphics[natwidth=7.78in, natheight=5.83in, width=5.5cm]{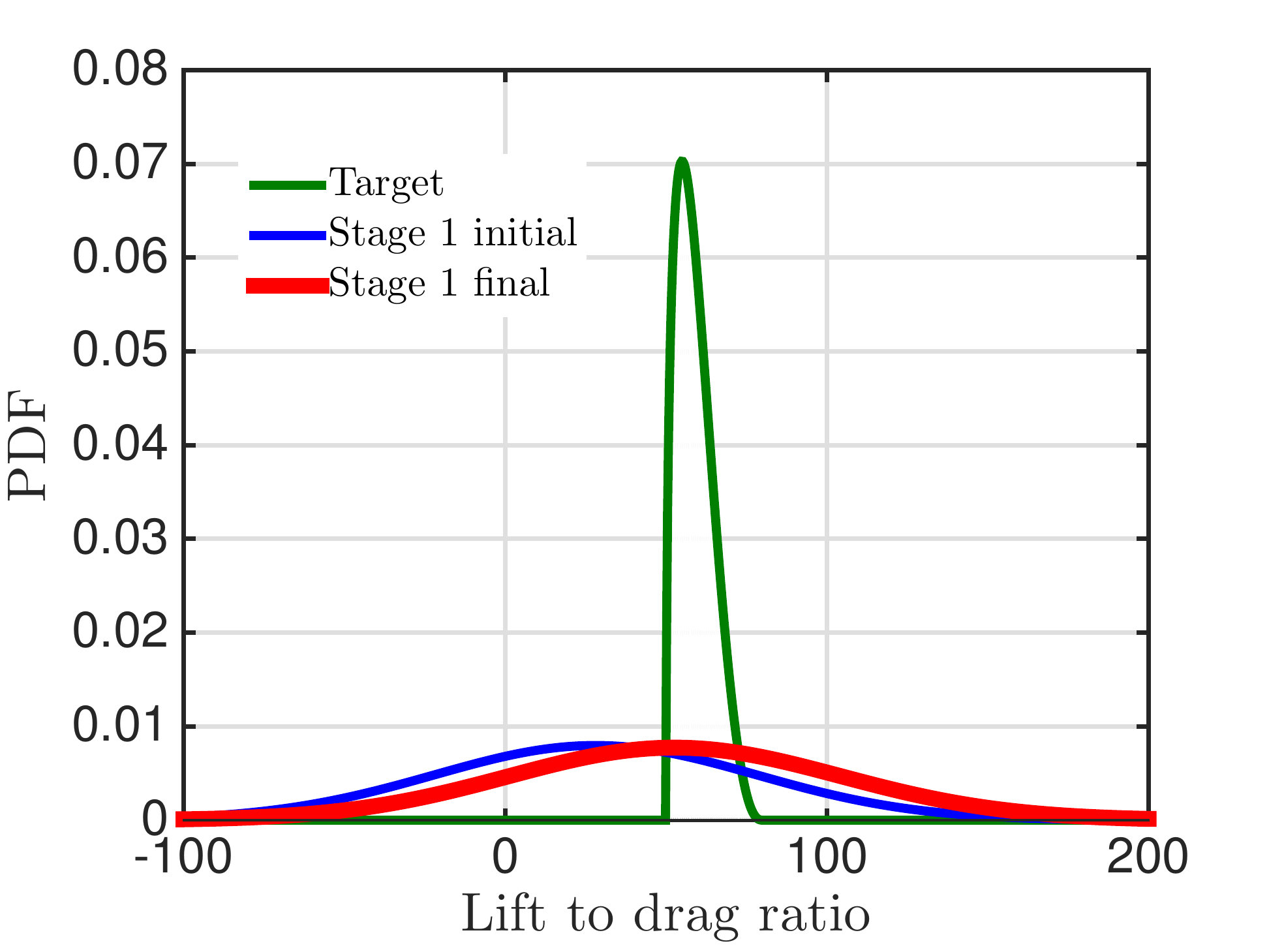}\label{fig:c1}}
\subfigure[]{\includegraphics[natwidth=7.78in, natheight=5.83in, width=5.5cm]{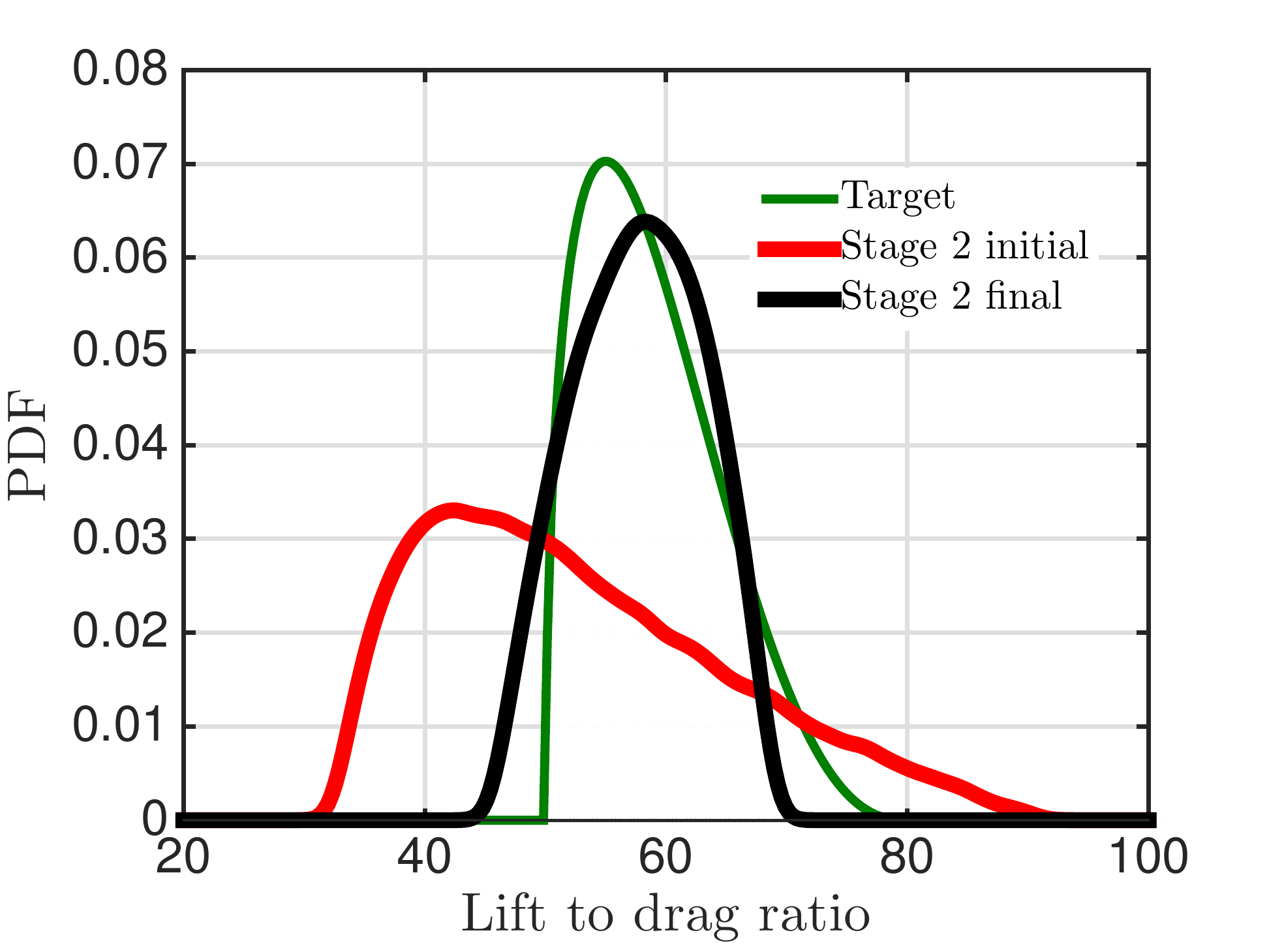}\label{fig:c2}}
\subfigure[]{\includegraphics[natwidth=7.78in, natheight=5.83in, width=5.5cm]{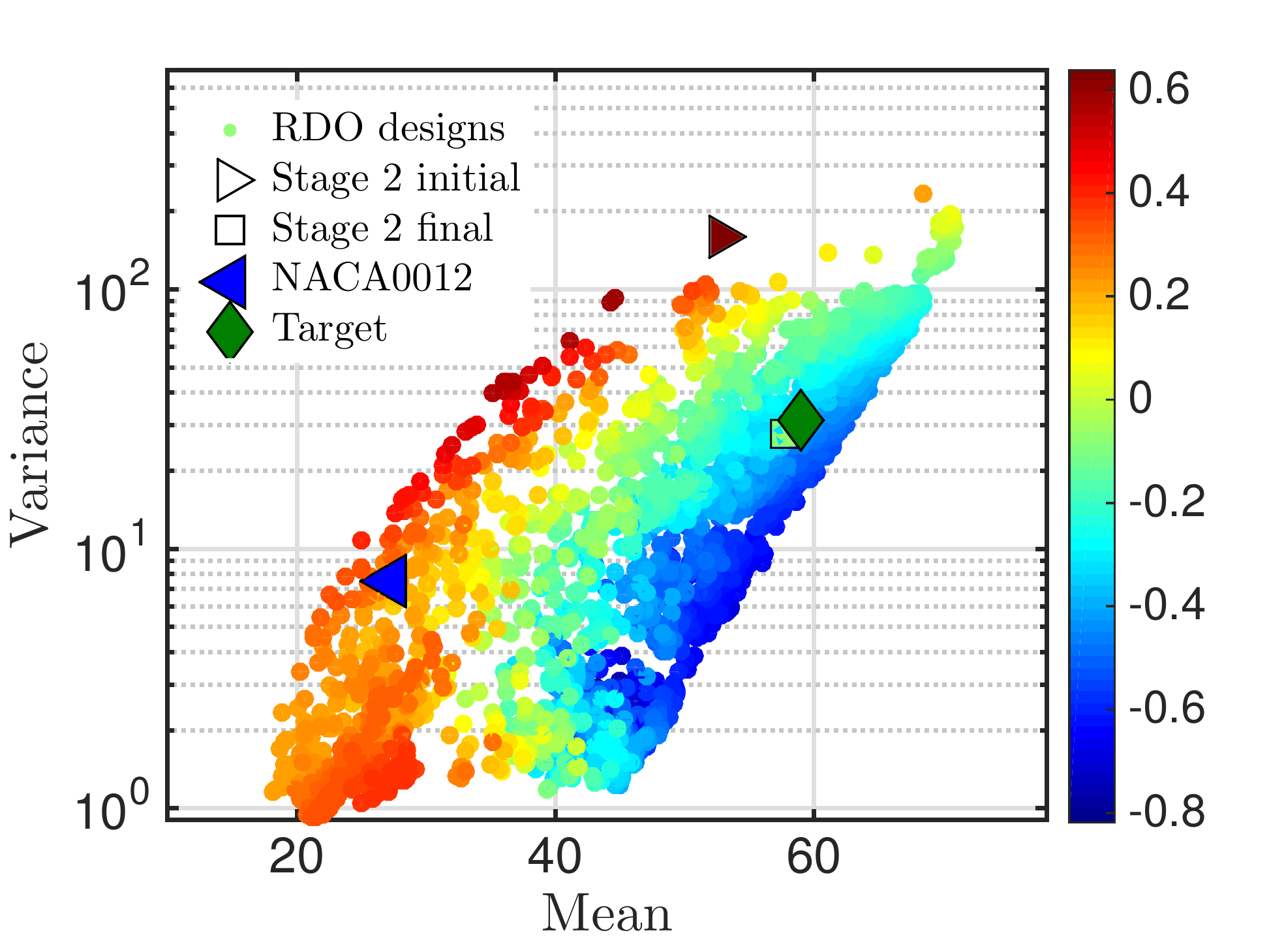}\label{fig:c3}}
\subfigure[]{\includegraphics[natwidth=7.78in, natheight=5.83in, width=5.5cm]{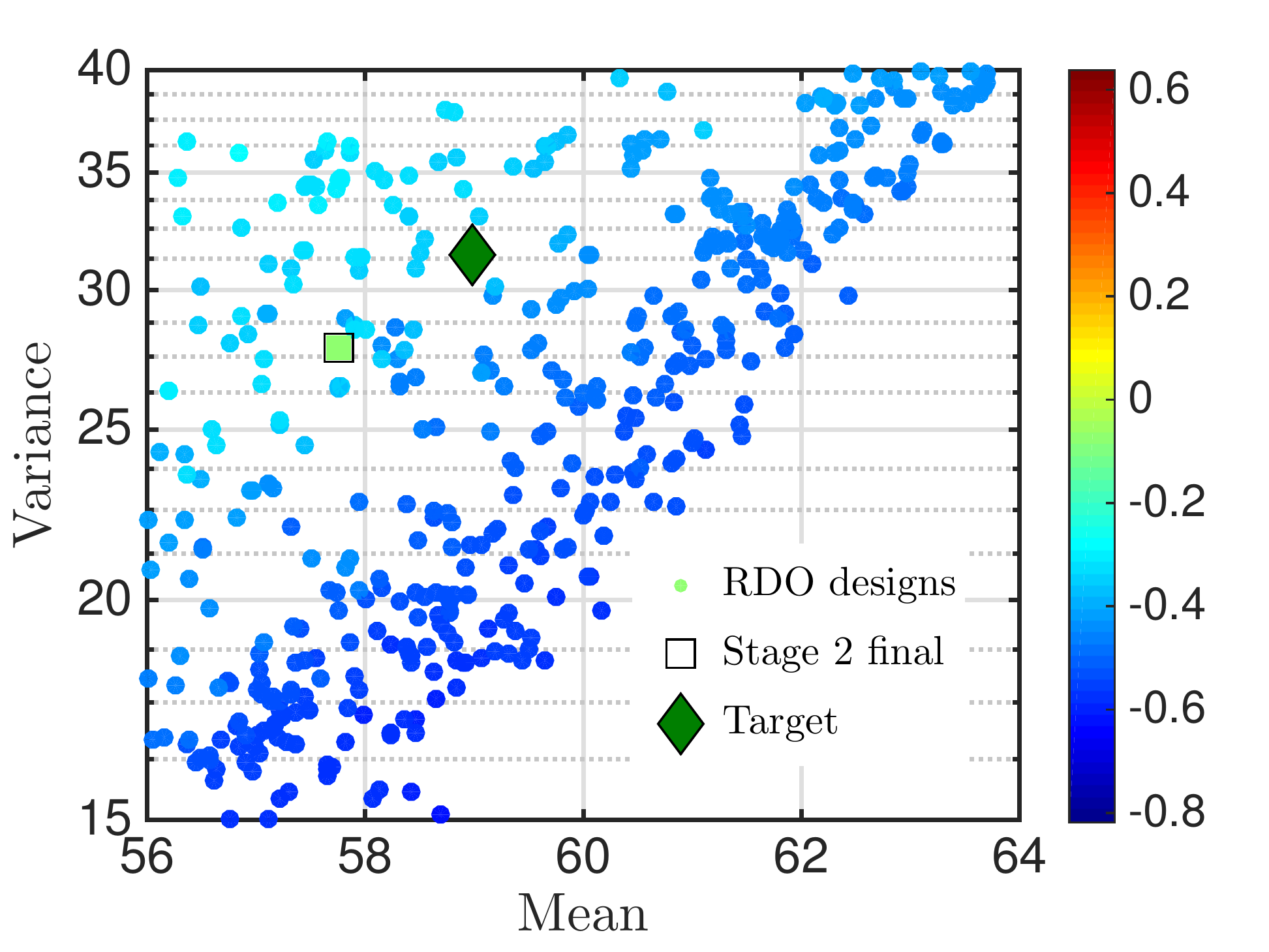}\label{fig:c4}}
\subfigure[]{\includegraphics[natwidth=7.78in, natheight=5.83in, width=5.5cm]{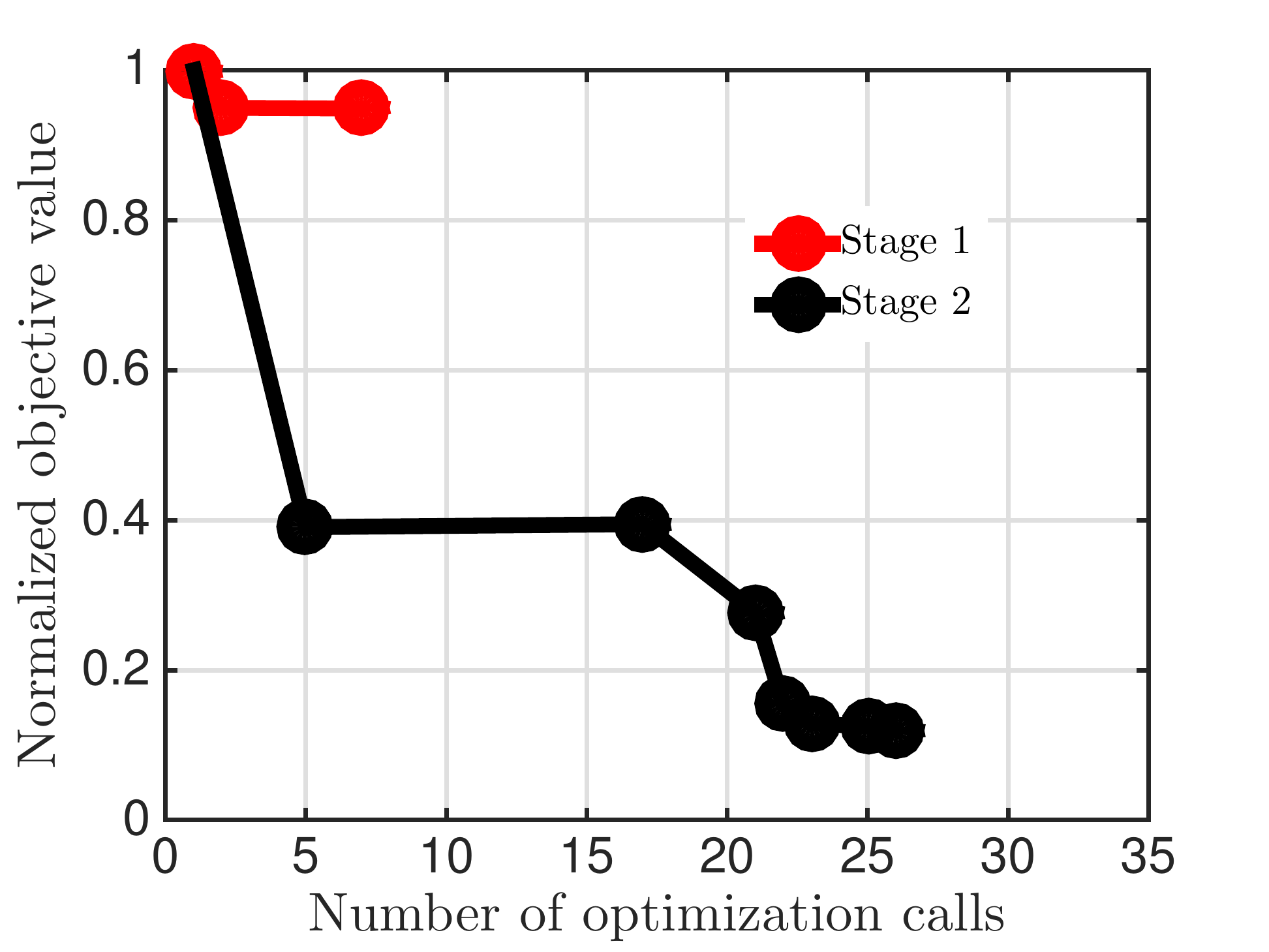}\label{fig:c5}}
\subfigure[]{\includegraphics[natwidth=7.78in, natheight=5.83in, width=5.5cm]{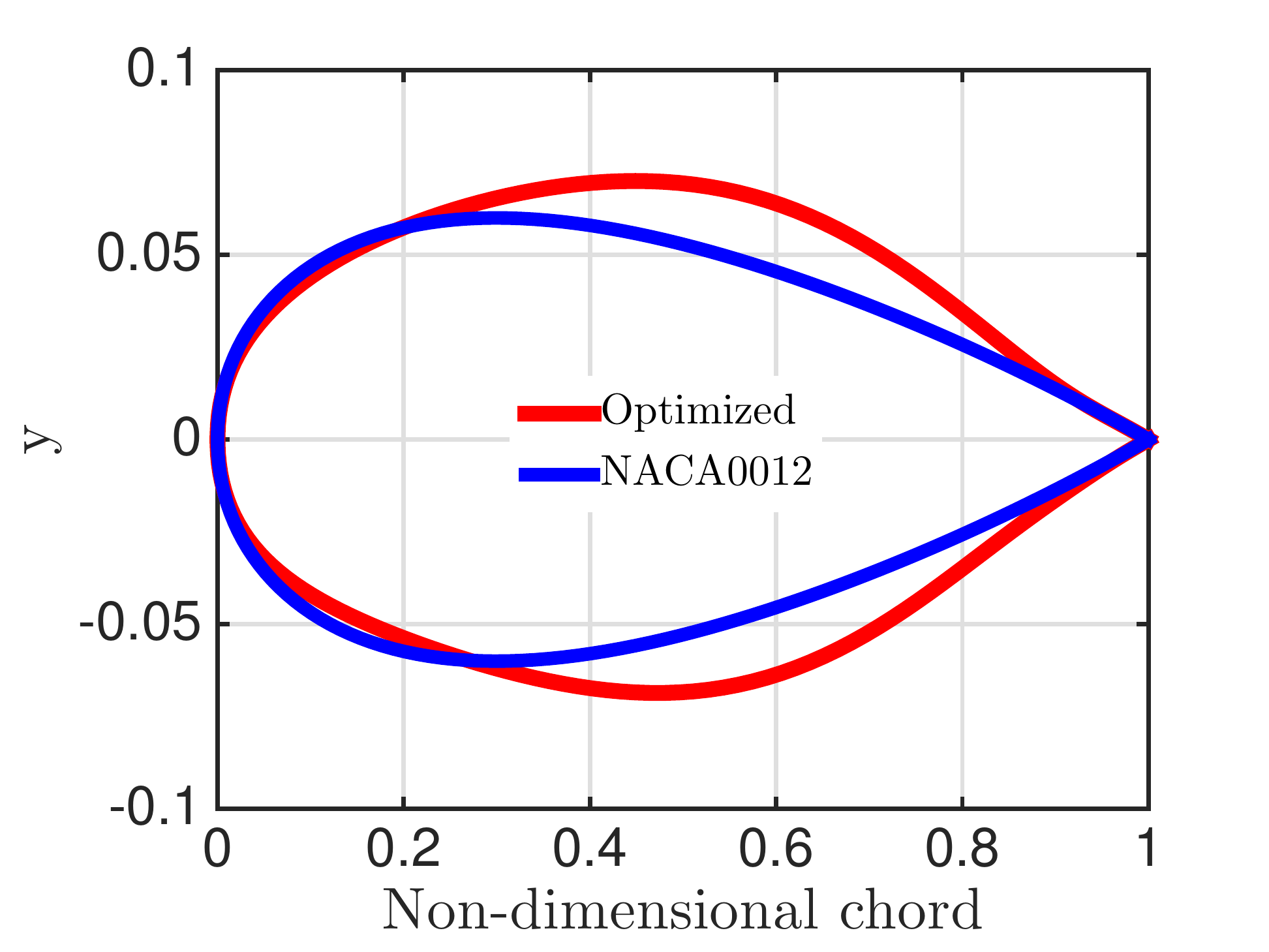}\label{fig:c6}}
\end{subfigmatrix}
\caption{Beta target results: (a) stage 1 ($h=50$), (b) stage 2 ($h=1$), (c) comparison with RDO designs with a close-up in (d), (e) stage 1 and stage 2 convergence plots, (f) stage 2 optimal design.}
\label{C}
\end{figure}

\noindent With beta density, we examine a case where the target is positively skewed. For this case, we use only one trial as in the Gaussian target case. The results are shown in Figure \ref{C}. The target has mean 59.0 and variance 31.5. The initial designs are shown in Figure \ref{fig:c1}, the final designs in Figure \ref{fig:c2}, and the convergence history in Figure \ref{fig:c5}. Stage 1 used 7 function calls with 3 major iterations, while stage 2 used 26 with 7 major iterations. For the beta target, stage 1 took 1 hour and 24 minutes while stage 2 took 5 hours and 8 minutes.

In the final result, we find that the optimizer tried to get as close as possible to the positively skewed target and produced a design whose density estimate has skewness -0.08. Neighboring solutions in the Pareto front---with similar means and variances---exhibited large negative skewness; see Figure \ref{fig:c3} and a closeup in Figure \ref{fig:c4}. This illustrates two key points. First, it demonstrates a limitation of using only the first two moments in RDO. At the same time, this highlights the strength of density matching. Using a simple $L_2$-norm distance metric, we push the optimizer to match not moments but the full pdf. This leads to a better match in all the moments of the response. The final design is shown in Figure \ref{fig:c6}. 

%% file: conclusions.tex
\section{Conclusions and Future Directions}
\noindent We present an alternative metric for optimization under uncertainty (OUU). We assume that the designer has provided a target pdf of the system response, and we minimize the distance between the design-dependent response pdf and the given target over possible designs. We study the differentiable $L_2$-norm between the response and target pdfs, though other distance metrics may be employed. One drawback of the $L_2$-norm is that if the target and response pdfs are not sufficiently large on the same support, then the objective function's gradient may not be useful for the OUU. We present a particular discretization of the $L_2$-norm objective function that uses a numerical integration rule and a kernel density estimate for the response pdf. The kernel density estimate with the Gaussian kernel has a simple form for the discretized objective's gradient. We offer two computational heuristics: (i) a two-stage strategy in the bandwidth choice for the kernel estimate that alleviates the support issue and (ii) a response surface approach for computationally expensive system responses. We apply this approach in two examples: (i) a simple function that produces a Gaussian response pdf whose variance is the design parameter and (ii) a CFD-based airfoil shape optimization. We show that the proposed pdf distance metric is a useful metric for OUU. 

Future work may explore the properties of different distance metrics for the response and target pdfs. We suspect that our study of the $L_2$-norm would translate easily to the Hellinger distance, which is a more common metric for comparing two pdfs. Also, it is possible to bound differences in moments by the Hellinger distance, which may enable more extensive and quantitative comparison with common RDO formulations. 

Another avenue for research may be to employ a maximum likelihood estimate on the sampled data to yield design pdfs with respect to an assumed distribution. Our process can still be applied in such a scenario. Moreover, in such cases the integrations can be performed analytically. An extension of this idea would be to use only Gaussian distributions, by estimating the mean and variance for all designs. This could atleast in theory lead to a more `classical' recovery of RDO for certain corner cases.